# Multi-Period Travelling Politician Problem: A Hybrid Metaheuristic Solution Method

**Masoud Shahmanzari**

*Faculty of Business,*
*Istanbul Commerce University,*
*İstanbul, Turkey*

**Deniz Aksen**

*College of Administrative Sciences and*
*Economics, Koç University,*
*İstanbul, Turkey*

**Saïd Salhi**

*Kent Business School,*
*Centre for Logistics and Heuristic Optimisation,*
*University of Kent, Canterbury, United Kingdom*

## ABSTRACT

This paper studies the *Multi-period Travelling Politician Problem* whose objective is to maximise the net benefit accrued by a party leader during a fixed campaign period. The problem is also characterised by flexible depots since the daily tours realised by the party leader may not start and end at the same city. A hybrid multi-start Iterated Local Search method complemented with a Variable Neighbourhood Descent is developed to solve the problem heuristically. Two constructive procedures are devised to generate initial feasible solutions. The proposed method is tested on 45 problem instances involving 81 cities and 12 towns in Turkey. Computational results show that the hybrid metaheuristic approach outperforms a recently proposed two-phase matheuristic by producing 7 optimal solutions and 17 new best solutions. In addition, interesting practical insights are provided using scenario analysis that could assist campaign planners in their strategic decisions.

**Keywords** Routing • Election logistics • Travelling politician problem • Iterated Local Search • Variable neighbourhood descent • Scenario analysis

## 1 Introduction

This paper looks into the *Multi-Period Travelling Politician Problem* (MPTPP) which arises in the planning of electoral campaigns. A set of cities with population- and time-dependent rewards and a set of campaign days are given. The politician's campaign team seeks to produce an efficient meeting schedule for the politician. Their objective is to collect the highest possible total reward from the meetings held in selected cities until the end of the campaign period. The schedule consists of either open or closed daily tours with the added restrictions of maximum tour duration and the necessity to return to the campaign centre frequently. The name of the MPTPP derives from *Periodic Travelling Politician Problem* originally coined in Aksen and Shahmanzari (2018). It can be classified as an application of the *Roaming Salesman Problem* (RSP) to election logistics. The RSP was first introduced in Shahmanzari et al. (2020) as a strategic campaign planning and reward collection problem with rich routing elements. In addition to capturing the MPTPP, the RSP can be used to model such problems as touristic trip planning, marketing campaign planning, and planning of client visits by company representatives.

The MPTPP which we investigate in this paper has some similarities with a multi-period extension of the *Prize Collecting Travelling Salesman Problem* (PCTSP) with time-dependent rewards and multiple visits. However, the tour does not have to be a round-trip. Besides, a city may be visited either in transit or for the purpose of a meeting.

There are two types of rewards, namely base rewards and depreciated rewards. The base reward is defined a priori for each city according to its characteristics, while reward depreciation circumvents successive meetings in the same city within short time intervals. It is also worth noting that cities of high population or significance may host multiple visits.

MPTPP has five unique features. First, a maximum tour duration is to be observed in each daily tour. Second, it is a selection problem since not every city has to host a meeting. Third, daily tours can be either an open or a closed route. Fourth, each city is associated with a time-dependent reward which changes linearly according to the day of the meeting and the recency of the previous meeting in that same city. Finally, certain cities may host more than one meeting.

Consider two sets of cities denoted by $\mathbf{V} = \{1,...,n\}$ and $\mathbf{N} = \mathbf{V} \cup \{0\}$ where the latter includes a fictitious city (indexed as 0) and the capital city (indexed as 1) which acts as the campaign centre. Also given is a set of days $\mathbf{T} = \{1,...,m\}$ before the elections. On each day $t \in \mathbf{T}$, any city $i \in \mathbf{N}$ can be visited with or without a meeting. A base reward $\pi_i > 0$ is specified for holding a meeting in city $i \in \mathbf{N}$. The travelling cost from city $i$ to city $j$ is known in advance and denoted by $c_{ij}$, $\forall i, j \in \mathbf{N}, i \neq j$. Similarly, the travelling time between each pair of cities is given by $d_{ij}, \forall i, j \in \mathbf{N}, i \neq j$.

There is a maximum allowed daily tour duration $\phi$. From a marketing and practical viewpoint, there is also an explicit limit on the number of meetings that can be held per day which is referred to as $\alpha$. There is a meeting time associated with each city $i \in \mathbf{N}$ which is denoted by $\mu_i$. The campaign of the politician is assumed to start in the base (capital) city $i = 1$ in the morning of day $t = 1$, and ends in the evening of day $t = m$. At the end of a day $t \in \mathbf{T}$, the politician stays overnight in some city $i \in \mathbf{N}$. We remark that waking up or staying overnight in a city $i$ does not necessarily yield a reward collection in that city. Also, for political reasons, the politician cannot be away from the capital for more than $\kappa$ consecutive days.

The MPTPP is a generalised version of the well-known travelling salesman problem (TSP) which is one of the oldest $\mathcal{NP}$-hard combinatorial optimisation problems. For several decades, the academic literature concentrated on developing heuristic approaches to find high-quality solutions for large-scale combinatorial optimisation problems. Recently, cross-fertilisation of different optimisation techniques—including not only (meta)heuristics, but also exact algorithms—has been introduced as a new variant of algorithms to deal with such problems. The motivation behind this effort is to exploit the potential synergy that might be harvested from the combination of complementary algorithmic characteristics. Especially the past 12 years saw a substantial growth in the numbers of workshops, conferences, books, special journal issues, and original research articles about hybrid algorithms the scope of which is not limited to combinatorial optimisation problems, but also extend to problems of continuous and multi-objective optimisation (Blum et al., 2011). The method we propose in this paper consists of two distinct metaheuristics to tackle the MPTPP, thus fits the category of hybridisation of such techniques among themselves.

As already mentioned above, Shahmanzari et al. (2020) recently solved the MPTPP in the framework of an election logistics application of the RSP by using a two-phase matheuristic. In this study, we propose the hybridisation of a multi-start iterated local search and a variable neighbourhood descent method as an improved solution technique for the same problem. The new solution approach produces better results. In addition, an extensive and methodical scenario analysis is conducted which returns a more in-depth perception of various situations that might arise in the planning of an election campaign. The contribution of this study is threefold.

1. Designing a hybridised metaheuristic solution method that integrates **I**terated **L**ocal **S**earch **(ILS)** and **V**ariable **N**eighbourhood **D**escent (**VND**) to solve the problem recently proposed by Shahmanzari et al. (2020)
2. Producing new best known solutions for large instances which can be used for benchmarking.
3. Presenting an extensive scenario analysis for obtaining managerial insights in the framework of a real-life election logistics application in Turkey.

The remainder of the paper is structured as follows. Section 2 reviews the literature relevant to the MPTPP. It is followed by Section 3 which describes the problem. Section 4 elaborates the proposed hybrid multi-start ILS algorithm and a VND metaheuristic collectively referred to as the MS-IVND algorithm. Detailed computational results are presented in Section 5 along with a scenario analysis in Section 6. The findings of the study are summarised in Section 7.

# 2 Literature review

The first TSP variant that is related to the MPTPP is the periodic travelling salesman problem (PTSP) where *n* cities have to be visited during a period of *m* days satisfying the required visit frequency of each city. One of the first formulations for the PTSP is provided in Cordeau et al. (1997). The objective is to minimise the travelling distance of the entire planning period while constructing a tour for each day and meeting the visit frequency of each city. Apart from additional assumptions in our model, the main distinction between PTSP and MPTPP lies in the fact that the rewards in MPTPP are time-dependent.

Moreover, three types of tours are present in MPTPP. Finally, unlike PTSP, MPTPP does not include a fixed depot.

Another class of TSP variants relevant to the MPTPP is known as *TSP with profits* (TSPP). This class of problems comprises *the prize collecting travelling salesman problem* (PCTSP), *the profitable tour problem* (PTP), and *the orienteering problem* (OP); see Feillet et al. (2005) for more information. MPTPP is more complex than these TSP variants as it contains multiple real life assumptions such as the presence of reward complexity, no penalty associated with non-visited cities, the requirement to visit campaign centre frequently, and not having a fixed depot. TSPP applications arise in a wide range of business operations, including realistic TSPs, job scheduling, freight transportation, or they occur indirectly as a subproblem in solution approaches dedicated to other routing problems. The TSPP is by definition the single criterion version of a bicriteria extension of the TSP with profit maximisation and travel cost minimisation being the two criteria. The basic characteristics of this generic problem are as follows.

(i). There is a value (like a profit or prize) associated with each vertex of the underlying graph.
(ii). A feasible solution is not required to visit all vertices.
(iii). A vertex can be visited at most once.
(iv). The distance (cost) matrix is nonnegative and satisfies the triangle inequality.

In PCTSP there are three attributes to consider, namely, the travel cost between cities *i* and *j*, the reward or the prize that is gained by visiting a given city and lastly the penalty of not being able to visit a given city. The aim is to determine a circuit that minimises the sum of travel costs and penalties of unvisited cities while guaranteeing at least a minimum total profit. Structural properties of the PCTSP related to the TSP polytope and the knapsack polytope were presented by Balas (1989) where families of facet-inducing inequalities were identified. Bounding procedures based on different relaxations were developed by Fischetti and Toth (1988) and Dell'Amico et al. (1995). The lower bound obtained according to the latter paper was used in a follow-up study by Dell'Amico et al. (1998) as the starting point of a Lagrangian heuristic. A branch-and-cut algorithm was proposed for the undirected PCTSP in the paper by Bérubé et al. (2009). The authors adapted and implemented some classical polyhedral results and derived cut inequalities for the PCTSP.

Another variant of the TSPP is the PTP where the objective is to maximize the net profit. This problem was first introduced by Dell'Amico et al. (1995). Fischetti et al. (2007) examined an extended, but related problem to cater for several vehicles with identical capacity. The last problem in the TSPP class is the OP. The goal in an OP is to determine a circuit or a path on a graph such that the sum of all collected prizes is maximised while still satisfying the upper bound on the total travel cost or time. Vansteenwegen et al. (2011) note that OP can be viewed as a combination of the *knapsack problem* and the TSP. Feillet et al. (2005) show that there is an equivalence between the path-seeking and circuit-seeking versions of the problem. The name *"orienteering problem"* originates from the treasure hunt game of orienteering in which individual competitors start at an initial control point, try to visit as many checkpoints as possible, and return eventually to the control point within a given time frame. Each checkpoint has its own reward. The objective of the game is to maximise the collected rewards. Several examples of OP applications have been cited in Ke et al. (2008) and in Vansteenwegen et al. (2011). These include orienteering competitions, routing technicians to service customers at geographically distributed locations, time-restricted fuel delivery to households with different urgency scores, athlete recruiting from high schools for a college team, pickup or delivery services with private fleets requiring the selection of only a subset of available customers, and trip planning for tourists visiting a city or a region. Another noteworthy OP application is found in Millar and Kiragu (1997). It involves fish scouting where a subset of fishing grounds are sampled to maximise the value of catch rate assessments. The authors referred to the underlying problem in this application as the *selective TSP* (STSP).

As mentioned in Section 1, MPTPP was first introduced in a paper by Aksen and Shahmanzari (2018) where the authors presented a mixed-integer linear programming formulation for the model solution of the problem. Later, Shahmanzari et al. (2020) proposed a two-phase matheuristic and newly generated data sets for a closely related version of the same problem, namely the RSP. We use the same data sets of the latter authors in our computational experiments. Aligned with the state-of-the-art literature, the variant that appears to be the most relevant and similar to our MPTPP is the multi-period OP with time windows which was investigated by Tricoire et al. (2010). The most critical difference is that in their study each tour starts and ends at the same central node known as the depot, whereas in our problem the terminal node of a tour is not known in advance. Also in our case, certain cities are allowed to be visited more than once, and the reward function changes according to the day of visit.

# 3 Problem description of MPTPP

The MPTPP is defined as follows. Given a set $\mathbf{N}$ of $n$ cities $(|\mathbf{N}| = n)$ and a set $\mathbf{T}$ of $m$ days $(|\mathbf{T}| = m)$, each city $i$ is associated with a nonnegative reward $\pi_i$ and a meeting duration $\sigma_i$. For each daily tour, there is a maximum duration $\phi$ and a maximum number of meetings $\alpha$. There are three possible types of tours the politician may consider:

*Type 1 tours: Multi-city closed tour*

The politician starts the day (wakes up) in city $i$ on day $t$, leaves $i$ and visits at least one more city scheduled for that day. At the end of the day, the politician returns to the same city $i$ to stay overnight. Type 1 tour is a closed tour starting and ending at city $i$ and involving at least one more city other than city $i$ (see Fig. 1).

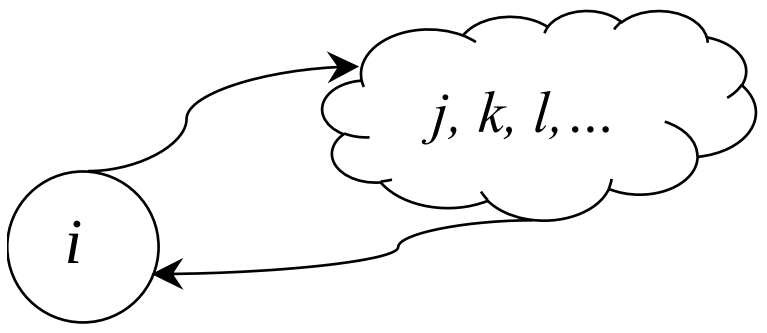

Fig. 1 Type 1 tour

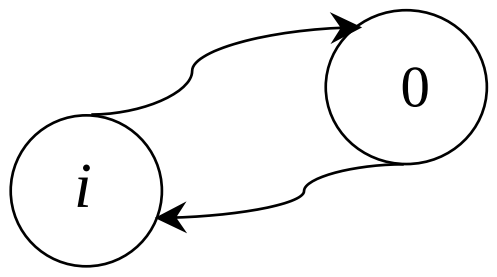

Fig. 2 Type 2 tour
(0: Fictitious City)

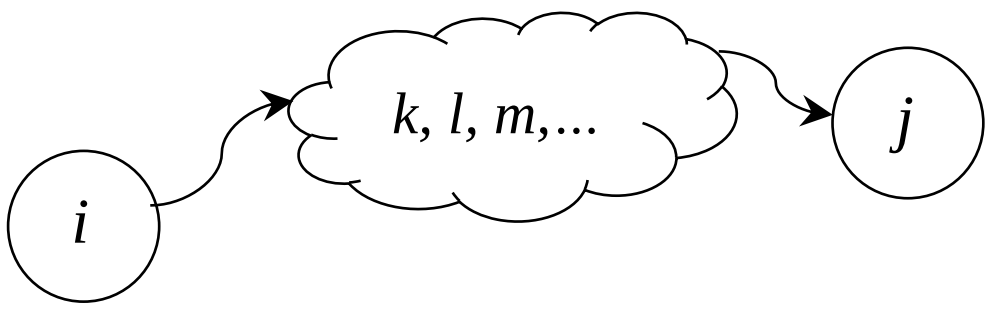

Fig. 3 Type 3 tour

*Type 2 tours: Single-city tour*

The politician wakes up in city $i$ on day $t$, spends the whole day in the same city collecting the meeting reward and stays overnight in the same city. In Type 2 tours we assume that the politician goes from city $i$ to the abovementioned fictitious city 0 and returns from 0 to $i$. This tour is therefore treated as a closed tour starting and ending at city $i$ (see Fig. 2).

*Type 3 tours: Multi-city open tour*

The politician wakes up in city $i$ on day $t$, and goes to another city $j$. In between cities $i$ and $j$ the politician may visit one or more cities, or may directly travel from $i$ to $j$ where the politician stays overnight. Type 3 tour is an open tour starting in city $i$ and ending in city $j$ as shown in Fig. 3.

*Objective function*

$$\max.\ NET\ BENEFIT = \sum_{i \in \mathbf{N}} \sum_{t \in \mathbf{T}} \pi_i \frac{m-t+1}{m} FM_{it} + \sum_{i \in \mathbf{N}} \sum_{t \in \mathbf{T}} \sum_{1 \le s < t} \pi_i \frac{m-t+1}{m} \times \frac{s}{Km} R_{its} - \bar{K} \sum_{i \in \mathbf{N}} \sum_{j \in \mathbf{N}} \sum_{t \in \mathbf{T}} c_{ij} X_{ijt} \qquad (1)$$

$X_{ijt}$, $FM_{it}$ and $R_{its}$ in (1) are binary decision variables which refer to traversing arc $(i, j)$ on day $t$, holding the first meeting of city $i$ on day $t$, and holding a meeting in city $i$ on day $t$ when the previous meeting in the same city was held $s$ days ago, respectively. The product resulting from the multiplication of the base reward $\pi_i$ with the binary variable $R_{its}$ in the second summation term of the objective function *NET BENEFIT* in (1) is depreciated further by a coefficient $K$. This depreciation is applied to successive meetings held in the same city. A simple normalisation coefficient denoted by $\bar{K}$ in (1) is used to make the rewards and travelling costs compatible. For completeness, the mathematical formulation of MPTPP is provided in Appendix A.

The reward that accrues from meetings is calculated according to the following rules. Further details are provided in Section 5.1 where we elaborate time-dependent rewards.

*Rule* 1: The earlier a meeting in the campaign period, the higher its reward.
*Rule* 2: The sooner a meeting is held after another in the same city, the more its reward is reduced.

The following constraints are considered in our implementation.

(i) *Chain feasibility:* Except for the last day, the terminal city of each day should be the same as the starting city of the next day.

(ii) *Maximum tour duration:* The total length of each daily tour should not exceed $\phi$.

(iii) *Return to campaign centre:* The campaign centre must be visited as a terminal node at least once every $\kappa$ days.

(iv) *The maximum number of meetings for different city categories:* There cannot be more than three meetings in big cities. This limitation reduces to two for regular cities. The remaining cities can host at most one meeting during the entire campaign.

(v) *No repeated meetings in the same city on the same day:* Each city can host at most one meeting every day. This applies only for closed tours.

(vi) *At most one meeting in big cities on the same day:* If there is a meeting in one of the big cities, there cannot be another meeting in the remaining big cities on the same day.

(vii) *Maximum number of meetings per day:* There cannot be more than $\alpha$ meetings each day.

# 4 Proposed metaheuristic

MPTPP can be formulated as a 0–1 integer program as stated in Shahmanzari et al. (2020) where a two-phase matheuristic is also proposed. Their approach is able to solve small instances of the problem in short CPU times. However, the overall solution quality deteriorated rapidly in medium and large size instances. In this study we develop a hybrid metaheuristic that incorporates ILS and VND to overcome the shortcoming of the matheuristic in Shahmanzari et al. (2020). We refer to the new hybrid method as the Multi-Start Iterated Variable Neighbourhood Descent (MS-IVND). The idea of multi-start optimisation proved to be promising in routing problems as explored by Reihaneh and Ghoniem (2018) for the distribution of pallets dispatched from a food bank to distant destinations of non-profit organisations.

## 4.1 General framework

We propose a metaheuristic based on ILS which calls a local search procedure iteratively by feeding a different starting point in an attempt to escape local minima. In brief, ILS has two main components, namely, *Perturbation* and *Local Search* (LS). Once an initial solution is constructed, the perturbation step diversifies the current solution by generating a new solution with a small modification to avoid the algorithm from being trapped in the same local minimum. At each iteration, a new solution is generated randomly by the perturbation mechanism which is then utilised by LS. See Lourenço et al. (2003) for a comprehensive tutorial of ILS, and Salhi (2017) for an overview of heuristic search including ILS.

The high-level architecture of ILS is sketched in Algorithm 1. The algorithm starts by generating an initial solution $S_0$. Inside the main loop of the algorithm, a perturbation mechanism and a local search are applied to this solution. If the resulting solution $S^{*\prime}$ satisfies the acceptance criterion, it replaces $S^*$. This procedure is repeated until the termination condition is met.

**Algorithm 1.** Basic ILS

1: $S_0$ =GenerateInitialSolution (See Section 4.4)
2: $S^*$ =LocalSearch( $S_0$ ) (See Section 4.6)
3: **Repeat**
4: $S'$ =Perturbation( $S^*$ ) (See Section 4.5)
5: $S^{*\prime}$ =LocalSearch( $S'$ )
6: $S^*$ =AcceptanceCriterion( $S^*$ , $S^{*\prime}$ )
7: **Until** the termination condition is met
8: **End**

## 4.2 The MS-IVND metaheuristic

Our proposed metaheuristic contains four additional features.

(i) The cyclic use of two novel constructive heuristics to build the initial solution from where the algorithm restarts after a certain number of iterations without improvement. Using this technique, we allow the algorithm to escape local minima.

(ii) The incorporation of a variable neighbourhood descent (VND) as a local search. See the review by Hansen and Mladenović (2003) for an overview of the VND method, its variants and applications.

(iii) Two diversification strategies: The first one starts the algorithm from a new initial solution that is constructed by either one of the two heuristics presented in Sections 4.4.1 and 0, respectively. The second one perturbs the best known solution to generate a new starting solution for the local search.

(iv) The incorporation of neighbourhood reduction to accelerate the search

The pseudocode of the MS-IVND method is given in Algorithm 2. The main loop of MS-IVND (lines 4–19) performs a sequential ILS. The initial solution of each ILS iteration is constructed using either the heuristic ESCC or SHRC which will be explained in Section 4.4. At every ILS iteration (lines 11–19) a perturbation mechanism is applied followed by a VND. The current ILS loop terminates when the maximum number of iterations without improvement $(\ell_{\max})$ is reached. Once the internal ILS loop stops, the algorithm generates a new initial solution with different characteristics to start over again. This way the entire process turns into a multi-start type method. The idea is to explore other parts of the solution space in order to avoid losing time in the VND step whenever it is trapped at a local optimum. The stopping criterion of the main approach is a maximum number of iterations denoted by $iter_{\max}$. Once this threshold is exceeded, both VND and MS-IVND loops will terminate, and the best found solution $S^*$ is returned. We describe the primary components of MS-IVND in the following sections. We first provide the solution representation, which is followed by the generation of the initial feasible solution, the perturbation mechanism, and the local search procedure. To accelerate the search we also resort to a granular neighbourhood reduction scheme which is used throughout the local search.

**Algorithm 2.** MS-IVND

$iter \leftarrow 1$
$\ell \leftarrow 1$
$S^* = \varnothing$ *// $S^*$ is initially empty set.*
**While** $iter < iter_{\max}$ **Do**
  **If** $iter$ is an odd number
    $S_{init} \leftarrow ESCC()$ (See Section 0)
  **Else**
    $S_{init} \leftarrow SHRC()$ (See Section 4.4.1)
  **End If**
  $S_{Best} \leftarrow VND(S_{init})$ (See Section 4.6)
  **While** $\ell < \ell_{\max}$ **Do**
    $S_{Pert} \leftarrow Perturb(S_{Best})$
    $S_{Temp} \leftarrow VND(S_{Pert})$
    **If** $NetBenefit(S_{Temp}) > NetBenefit(S_{Best})$
      $S_{Best} \leftarrow S_{Temp}$
      Append $S_{Best}$ to the set of $S^*$;
      $\ell \leftarrow 0$
    **End If**
    $\ell \leftarrow \ell + 1$
  **End While**
  $iter \leftarrow iter + 1$
**End While**
**Return** the maximum element of $S^*$

## 4.3 Solution representation

A solution of MPTPP is encoded as two 2-dimensional lists of nodes, one for routing schedules and one for meeting schedules. Each list consists of $m$ arrays (one for each day of the campaign period) where the sequence of visited cities is indicated. Fig. 4 represents an example of a campaign of $m = 3$ days with four cities being visited on the first day while only three meetings are held. The meeting for city 4 which is the last city visited on the first day is held on day 2. Therefore, that city becomes the starting city of the second day's tour. In the solution representation of MPTPP, the first list records the sequence of visited cities for each day regardless of whether there is a meeting there or not, whereas the second list records the meetings of each day. The order of nodes in the second list is not important since rewards earned from meetings are not dependent on the actual order of visits to cities within a day.

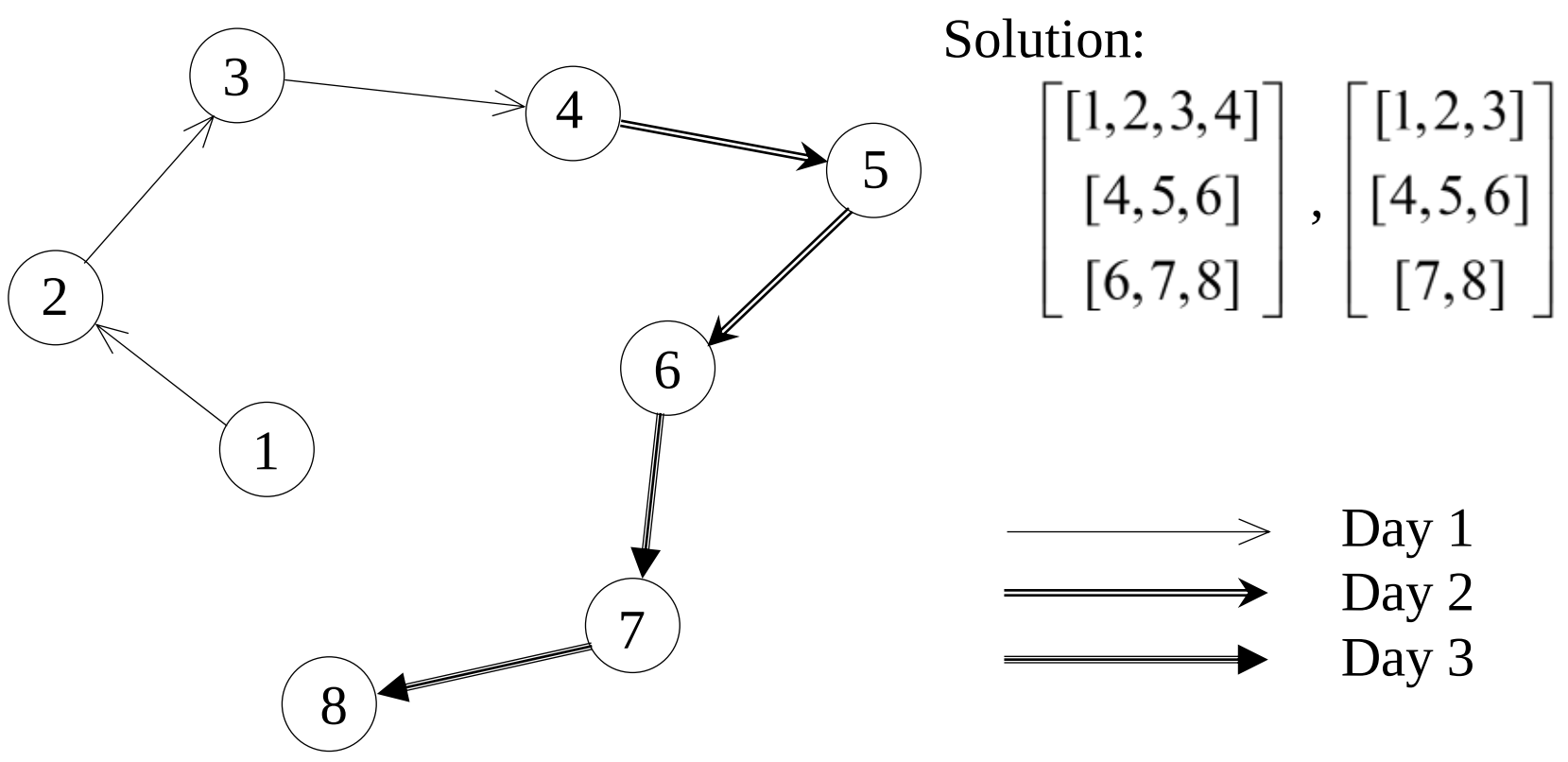

Fig. 4. Example of a solution with $n = 8$ and $m = 3$

## 4.4 Initial solution construction

We propose two constructive heuristics to obtain initial solutions for the MS-IVND. Every ILS in our implementation starts with a feasible initial solution that is produced by either a greedy heuristic or an exhaustive search. The main motivation for alternating between these two heuristics is to diversify the search. Both methods initiate with a solution consisting of *m* empty routes, each for one of the *m* days of the campaign period. Cities are iteratively inserted into routes until a stopping criterion is met. The first heuristic is called **S**electing **H**ighest-**R**eward **C**ities (**SHRC**) where it follows a similar schema to the nearest neighbourhood search. However, it inserts the next city into the route by considering the reward of the unrouted cities. The second heuristic which is named **E**xhaustive **S**earch of **C**andidate **C**ities (**ESCC**) performs an exhaustive search to find the best route for each day. We will describe these two procedures in the following two subsections.

### 4.4.1 The SHRC heuristic

SHRC starts by assigning a selection probability for each unrouted city; the higher the reward, the higher the chance to be selected for insertion. The first selected city is connected to the initial node, which is the campaign centre indexed as 1. After checking the feasibility of the route, the search either reverses to the previous move or jumps to another node. This procedure continues until the route's length exceeds the daily maximum tour duration. Then, the SHRC starts again to construct the route of the next day. In SHRC the highest-reward cities are mostly assigned to the early days of the campaign. The search continues until either no unrouted city can be inserted due to the violation of the maximum tour duration or no other city is available for insertion. The main steps of SHRC are presented in Algorithm 3. We assume that every city is visited only once, and every visited city hosts a meeting. For the terminal (depot) node of each day, it is possible to hold the meeting at the end (beginning) of that day or shift the meeting to the beginning (end) of the next (previous) day.

**Algorithm 3.** The construction of the initial feasible solution: **SHRC**

**For** $day = 1$ to $m$ **Do**
    **While** $S_{init}[day]$ is time feasible **Do**
        $i \leftarrow$ Pseudo-random selection of an unrouted city
        Append $i$ to $S_{init}[day]$;
    **End While**
    Drop $i$ from $S_{init}[day]$;
    $S_{init}[day+1] \leftarrow i$
**End For**
**Return** $S_{init}$

### 4.4.2 The ESCC heuristic

ESCC is also a greedy heuristic which builds daily routes one by one; but it seeks to assign cities to days relatively quickly. The detailed pseudo code of ESCC is presented in Algorithm 4. The main idea is to assign the cities with higher rewards to the early days of the campaign due to the characteristics of the reward function. This assignment procedure is performed through an exhaustive search repeated for each day of the campaign period (the outermost **For** loop in lines 2–28). For each day $t \in \mathbf{T}$, we first create a sorted list of all cities with respect to their updated rewards (lines 3–9). Then, we build a feasible route such that the maximum possible net benefit is achieved considering limitations like the maximum tour duration and the necessity to return to the campaign centre periodically (lines 10–23). To this end, all possible permutations of selected cities are evaluated within a CPU time limit of $MaxTimeInit$ for each day (lines 20–22). In our experiments we set $MaxTimeInit$ to one second. Finally, a feasibility restoration function is utilised to make the solution feasible (lines 24–26).

In the ESCC heuristic we assume that there will be a meeting in every visited city of the daily tour other than the wakeup city. It is interesting to note that the wakeup city of day $t$ never hosts a meeting; but being the terminal city for day $(t-1)$ where the politician stays overnight, it always hosts a meeting on day $(t-1)$.

**Algorithm 4. The** construction of the initial feasible solution: **ESCC**

Initialize WakeupCity(1) $\leftarrow$ *City* 1 and TopCities $= \varnothing$; // *City* 1 *refers to the campaign centre*.
**For** $t = 1$ to $m$ **Do**
  **If** $t = 1$
    Rewards $\leftarrow \pi$
  **Else**
    Calculate the reward of each city by taking into account the current meeting day $t$ and the recency of the previous meetings which may have been held before day $t$;
  **End If**
  Sort all cities eligible for hosting a meeting on day $t$ in descending order of their rewards into the array $\mathbf{H}_t$;
  Append WakeupCity($t$) to the set TopCities;
  $k \leftarrow 1$
  **While** $|\text{TopCities}| < 6$ **Do** // *We allow at most six visits on a tour.*
    $j \leftarrow \mathbf{H}_t[k]$
    **If** WakeupCity($t$) $= j$
      $k \leftarrow k + 1$
    **Else**
      TopCities $\leftarrow$ TopCities $\cup \{j\}$
      $k \leftarrow k + 1$
    **End If**
  **End While**
  **While** CPU time elapsed < *MaxTimeInit* **Do**
    Evaluate all the 5! tour permutations and select the one with the lowest total travelling cost. Break ties arbitrary; // *Each possible permutation represents an open tour for day t.*
  **End While**
  Call the selected tour BestTour($t$) $= [1^*, 2^*, 3^*, 4^*, 5^*, 6^*]$ // $1^*$ *is going to be WakeupCity*($t$).
  **If** BestTour($t$) is infeasible with respect to *MaximumTourDuration*
    Crop BestTour($t$) from its right end starting at city $6^*$ until it becomes time-feasible;
  **End If**
  WakeupCity($t+1$) $\leftarrow$ last visited city in BestTour($t$)
**End For**
**Return** BestTour($t$), $t = 1, ..., m$.

## 4.5 The perturbation mechanism

Perturbation plays a crucial role in our implementation since it forces the algorithm to escape from the present local optimum. One of the main goals in such schemes is to control the level of perturbation by avoiding too strong or too loose modifications as this balance helps to maintain the desirable properties of the current solution. The perturbation phase used in MS-IVND accepts both improving and non-improving moves. In our implementation, the perturbation operator returns a randomly perturbed solution by considering three potentially-deteriorating moves which are randomly selected in each iteration.

M1: Two pairs of two distinct cities are randomly selected and their positions are swapped.

M2: The two cities with the lowest rewards are replaced with two unrouted cities.

M3: A city from a route is randomly removed and inserted into the cheapest position in another route.

It is worth noting that the above perturbation transforms the current solution configuration into a new configuration that is unlikely to be reversed to its previous state. This sort of a perturbation mechanism is commonly used in tabu search to avoid cycling. We have also tested other moves such as the well-known 2-Opt. In this classical move, two edges in a route are removed and reordered to eliminate crisscrosses which would augment the total travelling distance when triangular inequality applies (Croes, 1958). However, the results were found to be inferior; and hence not used here.

## 4.6 Local search procedure

We use the variable neighbourhood descent (VND) method as our local search procedure (LS). VND searches for improving solutions within various defined neighbourhoods used in a systematic way to avoid local optimality (Hansen & Mladenović, 2003). VND is in principle similar to, though much simpler than, the multilevel metaheuristic developed by Salhi and Sari (1997) for multi-depot routing problems where local searches known as levels are used instead of neighbourhoods. The main steps of VND are sketched in Algorithm 5.

**Algorithm 5.** VND

$Improved \leftarrow \text{FALSE}$
**While** $Improved = \text{FALSE}$ **Do**
    $k \leftarrow 1$
    **While** $k \le k_{\max}$ **Do**
        $k' = ((k-1) \bmod 4) + 1$
        $S^* \leftarrow LocalSearch(S, N_{k'})$
        **If** $NetBenefit(S^*) > NetBenefit(S)$
            $S \leftarrow S^*$
            $Improved \leftarrow \text{TRUE}$
        **End If**
        $k \leftarrow k + 1$
    **End While**
**End While**
**Return** $S$

### 4.6.1 Neighbourhood structures

We introduce four neighbourhood structures with two being inter-route and two intra-route.

*$N_1$:* 1-1 *Exchange Inter-route*:

Two nodes from two distinct routes are randomly selected, and their positions are exchanged.

*$N_2$:* *Drop-Add*:

One node is randomly selected and dropped from its current route. Next, the unrouted node with the highest reward is selected and inserted into the cheapest position in other routes.

*$N_3$:* 1-1 *Exchange Unrouted*:

Two nodes are randomly selected from the set of routed and unrouted nodes, respectively. The interchange is made between the two.

*$N_4$:* 1-1 *Exchange Intra-route*:

Two nodes in the same route are randomly selected and their positions are swapped. The selected nodes can be either a depot (terminal) node or a transient city. In case of a depot (or a terminal node) is selected, the last (first) node of the previous (next) period is changed as well.

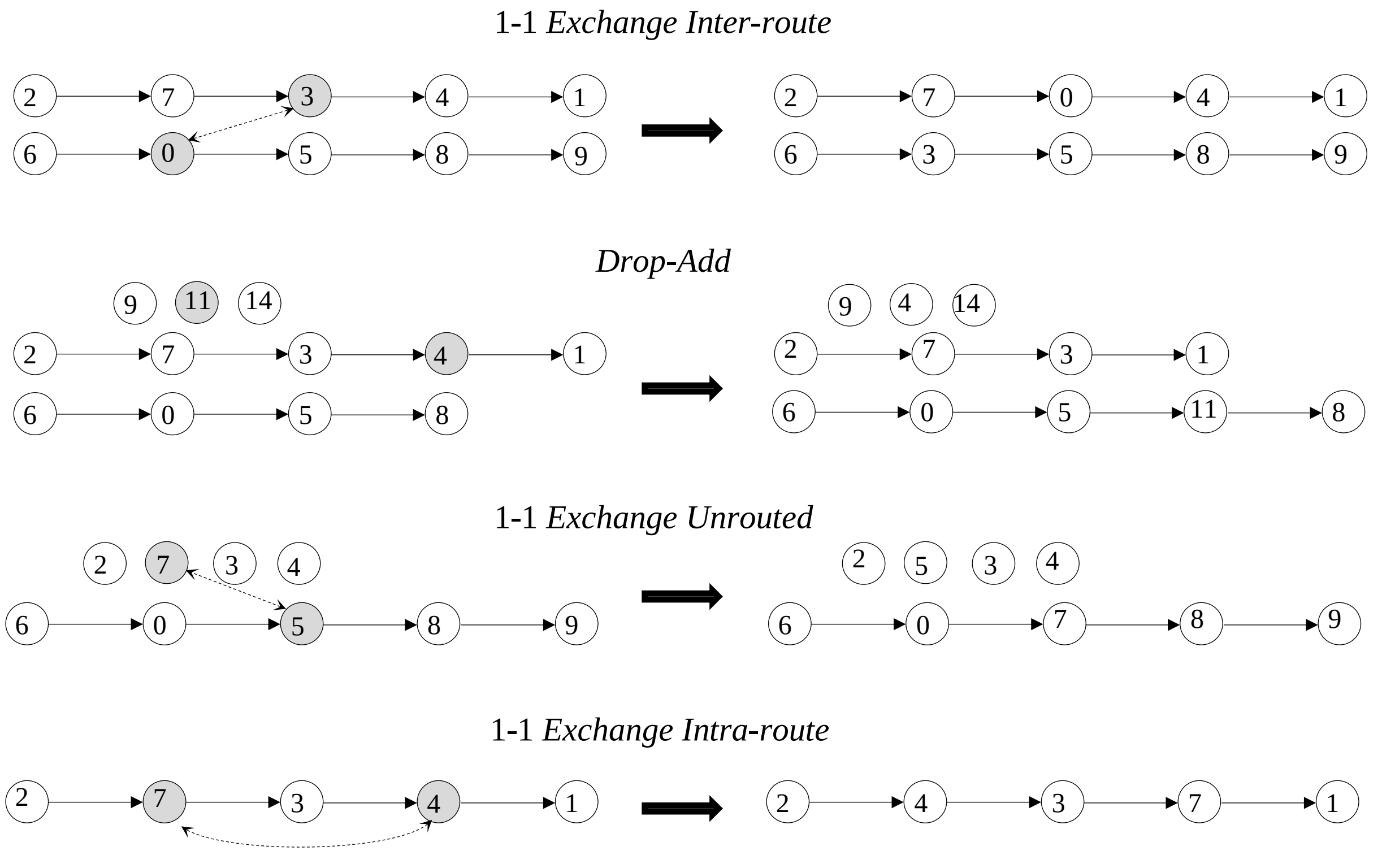

Fig. 5. Examples for the four moves used in VND

Fig. 5 depicts these four moves. As for all neighbourhoods, we execute the best improving moves (if any) where they might affect either scheduling or routing decisions. Each operator is able to select any type of nodes including depot, terminal, and transient nodes. For example, in the first two cases, the moves are more complicated due to the resulting changes in the route of the next or the previous period. Moreover, if a node is inserted as a depot (terminal node), we check whether holding the meeting in the current day or the previous (next) day improves the objective value. If this is the case, we assign the meeting to that day.

In our hybrid MS-IVND method, we explore only feasible regions of the solution space since the feasibility restoration is quite time-consuming for the MPTPP. Hence, a solution is only accepted if it is feasible with respect to the constraints defined in Section 3.

### 4.7 Acceleration scheme

For the purpose of accelerating the search, we incorporate a neighbourhood reduction scheme similar to granular neighbourhoods in Tabu Search (Toth & Vigo, 2003). During the insertion of node $i$ into route $R(j, j+1, \ldots, k-1, k)$, if we observe $\min(t_{ij}, t_{ik}) \leq \lambda t_{jk}$ where $\lambda$ is an appropriate positive coefficient and $t_{ij}$ denotes the travel time between $i$ and $j$, we allow the evaluation of this move. Otherwise, this move is eliminated. The logic behind such a restriction is to prevent the algorithm from consuming CPU time by evaluating non-promising moves. Obviously, if node $i$ is too far from the starting and ending nodes of a route, such an insertion becomes a non-promising move. Therefore, in our MS-IVND implementation only those moves are taken into consideration in which $i$ is fairly close to either $j$ or $k$ or both. Similar neighbourhood reduction schemes were successfully developed earlier for the multi-depot vehicle routing problem by Salhi and Sari (1997). These acceleration schemes are also highlighted in Salhi (2017) as part of an effective design in heuristic search in general, and fully exploited in Sze et al. (2016, 2017) for the classical and the cumulative vehicle routing problems.

## 5 Computational results

We tested our hybrid solution method MS-IVND on a Dell Precision T7810 model PC equipped with Intel Xeon® E5-2690 v4 2.60 GHz processor and 32 GB RAM running 64-bit Windows 10 Professional. The Xeon® processor provides 28 threads with the hyper-threading feature turned on. We coded the mathematical models of MPTPP and the MS-IVND method in Python 3.6.3 (64-bit version). For the solution of the models we employed the commercial mixed-integer linear programming (MILP) solver GUROBI 7.5.0 (GUROBI Optimization, 2020) which is called within the Python code.

The solver specific options applied to all runs are listed in Appendix B where more information is provided. Using 45 MPTPP instances from Shahmanzari et al. (2020) we compared the performance of MS-IVND with the best known solutions and commercial solver results in that paper. These instances were generated with real-world travel distances and times among all cities of Turkey. In each instance at least one big city, one midsize city and one small city are included. Appendix C presents the characteristics of a medium-size instance with $n = 40$ cities and a planning horizon of $m = 10$ days. The best tours for this instance are mapped in Fig. C.2. All instances are publicly available at http://shahmanzar.ir/MPTPP.html.

We consider the 80 cities (provinces) of Turkey plus a campaign centre, namely, the capital city Ankara. At the time of the 2015 elections, Turkey had 81 cities and 85 electoral zones where İstanbul was comprised of three zones, İzmir and Ankara of two zones each. Each city is associated with a base reward value and a fixed meeting duration.

The maximum tour duration limit $\phi$ which is 14 hours (12 hours) in the summer (in the winter) imposes an implicit threshold on the number of cities that can be visited any given day. Each city can accommodate at most one meeting a day. There can be an upper bound (such as two or three) on the total number of meetings held in each city during the campaign period. The meeting durations range from 60 to 120 minutes depending on the population of the host city. For the three biggest cities, namely İstanbul, Ankara and İzmir, it is 120 minutes. For cities with fewer than one million population it is 60 minutes, and 90 minutes for all other ones. Another point to be mentioned is the periodic returns to the campaign centre Ankara. The politician cannot be away from the capital city for more than $\kappa$ consecutive days with $\kappa \in \{4,5,6,7\}$.

The travel time between a pair of cities is the fastest trip (either by road or by airplane), and its cost is the sum of the monetary costs of the trip legs. The reader interested in the derivation of these two entities is referred to Shahmanzari (2019) where the trade-off between road and air travel is explained in detail.

## 5.1 Time-dependent rewards

In this section we expound the reward calculation and the categorisation of cities in Turkey from the main opposition party's perspective. The proposed model utilises a multifaceted reward function. Initially, a nonnegative prize of $\pi_i$ (*base reward*) is specified for holding a meeting in each city $i \in \mathbf{V}$ where $\pi_i$ depends on two factors:

Factor 1: Population of city $i \in \mathbf{V}$ $(Pop_i)$.

Factor 2: Ratio of votes of the politician's party (PP).

In addition, the reward earned in city $i \in \mathbf{N}$ on day $t \in \mathbf{T}$ is dependent on two further factors:

Factor 3: Number of remaining days denoted by $(m-t)$ until the end of the campaign.

Factor 4: Number of days passed since the previous meeting in the same city, denoted by $s$ where $1 \le s \le t-1$.

The first two factors directly affect the base reward $\pi_i$, while the remaining two make the reward time-dependent. Each factor is explained below.

Factor 1: Population

In an ideal representative election system, the number of representatives allocated to each electoral zone (to each *city* in the context of our problem) has to be approximately commensurate with the population share of that zone (Çarkoğlu & Aksen, 2019). This implies that population is one of the most decisive factors in determining the importance of a city in an election campaign.

In the calculation of the parameter $\pi_i$, each city is first assigned a starting base reward of 100 units. Each city's population is divided by the minimum population of all cities, rounded to the nearest integer, and multiplied by a city-dependent multiplier. The result is added to the initial reward of 100, and then multiplied further by a *Criticality Factor* $(CF_i)$. The formula of the base reward calculation is shown in (2) where where the operator $[\![\cdot]\!]$ rounds its argument to the nearest integer number.

$$\pi_i = CF_i \times \left(100 + \left[\!\!\left[ \frac{Pop_i}{Min.Pop.} \right]\!\!\right] \times Multiplier_i \right) \tag{2}$$

The multiplier in (2) is taken as 3.0 for the top seven (most populated cities), but 2.0 for İstanbul. The remaining cities are assigned a multiplier value of 5.0. When a town is taken into consideration, it is assumed that it shares the same multiplier with its parent city. These specific multipliers turned out to yield proportionate and scaled base rewards.

Factor 2: Ratio of votes

In order to find the effect of variation in the number of votes on the number of deputies in the parliament, the data of June 2015 election has been analysed for all cities. In our criticality analysis, we first simulated the election procedure according to the actual vote counts registered in the election of June 2015. We were able to reproduce exactly the same seat distributions in all 85 electoral zones of Turkey which shows the validity of the implemented simulation. Next, we evaluated each city by decreasing and increasing the votes of the political party in that city by 20%. Cities are categorised as discussed below. The reward statistics are provided in Table 1.

Table 1 Statistics of rewards in criticality categories

| | Noncritical | Negative Critical | Positive Critical | Positive-and-Negative critical |
|---|---|---|---|---|
| Number of cities | 42 | 19 | 11 | 9 |
| Average Reward | 268 | 444 | 564 | 1,193 |
| Min. Reward | 210 | 315 | 460 | 800 |
| Max Reward | 440 | 675 | 680 | 2,370 |

We define four criticality categories to label the importance of a city from the perspective of the political party on the basis of the distribution of votes in the previous election. Different towns or electoral zones of a city are mutually assigned to the same criticality category. We define four distinct categories of criticality of cities and associate each category with a different Criticality Factor $CF$ as follows:

*Category 1*: *Noncritical Cities*

These are the cities in which the number of seats won by the political party would not change even when the number of its votes changes by ±20%. We set $CF_i = 2$ for $i \in Noncritical\ Cities$.

*Category 2*: *Negative Critical Cities*

In these cities, a 20% increase in the votes of the political party does not affect its seat number in the parliament (the number of its deputies elected from those cities). However, a 20% decrease would cause the political party to lose at least one seat. We set $CF_i = 3$ for $i \in Negative\ Critical\ Cities$.

*Category 3*: *Positive Critical Cities*

In positive critical cities, the political party would gain at least one more seat in the event of a 20% increase in the vote count of the past election. However, there exists no risk of losing any seat in the event of a 20% loss in the votes. We set $CF_i = 4$ for $i \in Positive\ Critical\ Cities$.

*Category 4*: *Positive-and-Negative Critical Cities*

The situation is most sensitive in cities of category 4 where an increase or decrease by 20% in the vote count would impact the party's current seat count in the parliament. Hence, we set $CF_i = 5$ for $i \in Positive\text{-}and\text{-}Negative\ Critical\ Cities$.

The motivation behind these *CF* values is to assign high rewards to highly populated cities, but doing so at a decreasing rate. Another motivation is to close the enormous gap between metropolitan cities and other midsize cities of Turkey. For instance, İstanbul, despite its ~15 million population, should not earn thrice as much base reward as Ankara just because of having thrice as much population.

*Criticality Analysis*

To illustrate the effect of the *CF*, let us consider two cities, namely Samsun and Kastamonu in the Black Sea Region. The base reward of Samsun is higher than Kastamonu ($\pi_{Samsun} = 540$ and $\pi_{Kastamonu} = 500$) although Samsun's population is three times higher than Kastamonu's population. The base reward of Kastamonu almost catches up with Samsun because the latter is a positive-critical city, whereas the former is a negative-critical city. As highlighted earlier, MPTPP has obviously a selective nature where not all cities in $\mathbf{N} = \{0,1,...,n\}$ need to be included in the meeting schedule.

Factor 3: Number of remaining days until the election day

We assume that as we get closer to the end of the campaign, the influence of meetings will decrease. In order to inflate the base rewards with the increasing number of remaining days until the elections, we develop the following reward function.

$$Reward_i(t) = \pi_i \frac{m-t+1}{m} \qquad i \in \mathbf{N},\ t \in \mathbf{T} \tag{3}$$

If the political party decides to reverse the effect of Factor 3, the formula in (3) can be easily modified by setting $Reward_i(t) = \pi_i \frac{t}{m}$. The reward of a meeting would then be the lowest on the first day and the highest on the last day of the campaign.

Factor 4: Number of days passed since the previous meeting

In order to prevent the model from visiting highly rewarded cities frequently, we severely penalise repeated meetings. To inflate the base rewards with the increasing number of days passed since the last meeting, we extend the reward function previously given in (2) as follows:

$$Reward_i(t,s) = \pi_i \frac{m-t+1}{m} \times \frac{s}{Km} \qquad i \in \mathbf{N},\ t \in \mathbf{T} \tag{4}$$

where $s$ represents the number of days passed since the last meeting and $K$ $(K \geq 1)$ is a prespecified depreciation factor for repeated meetings. The criterion of depreciation is not the number of meetings held in city $i$ so far, but the recency of the previous meeting.

## 5.2 Parameter calibration

The efficiency of most metaheuristics depends on their corresponding parameters. The strength of our MS-IVND lies in the fact that it has three parameters only: (i) Maximum number of iterations $(iter_{\max})$, (ii) Maximum number of iterations without improvement $(\ell_{\max})$, and (iii) the granularity coefficient $(\lambda)$. As we use the same neighbourhood structures throughout the experimentation, the number of neighbourhood structures, namely $k_{\max}$ is fixed and hence not considered in our parameter calibration.

Starting with a promising configuration, we performed empirical tests to determine the best tuning parameters for MS-IVND on a benchmark set of 10 instances. Parameters were varied one at a time before the method was run again to observe the effect of a given parameter. We compared the test results obtained before and after varying the value of the parameter, and chose the one which yielded better results. If no change was suggested at the end of the tests conducted for a particular parameter, we moved to the new configuration. The following values have been tested: $iter_{\max} \in \{25, 50, 100, 150\}$, $\ell_{\max} \in \{25, 50, 75, 100\}$, and $\lambda \in \{0.25, 0.40, 0.75, 0.90\}$.

After the calibration tests we used the following parameter values in the sequel of our experiments: We assigned a fixed value of 100 to $iter_{\max}$ since the number of calls to the algorithms ESCC, HRC and VND is a significant determinant of the execution time of MS-IVND. The parameter $\ell_{\max}$ is set to 50 since the multi-start nature within the scope of MS-IVND allows the search to explore diverse regions of the solution space. Finally, the parameter $\lambda$ is fixed at 0.75. Comprehensive results of the calibration tests are not provided here, but can be collected from the authors. Due to the inherent randomness in our algorithm, the reported MS-IVND solution values and running times for every instance of MPTPP are the averages of the best results over ten runs.

## 5.3 Results on existing instances

To our knowledge, no metaheuristic has been published for the MPTPP. The only solution methodology developed is FDOR (**F**inding **D**aily **O**ptimal **R**outes) which is a matheuristic due to Shahmanzari et al. (2020). Hence, we benchmarked MS-IVND against the commercial MILP solver GUROBI and the FDOR heuristic on the same test instances as those reported in Shahmanzari et al. (2020). They are divided into three sets referred to as `PE.I`, `PE.II` and `LE`, which comprise 22, 20 and 3 instances, respectively. Our results are reported in Table 2 through Table 5. Instance names in the first column of each table comply with the identification template $n$`C-`$m$`D` where $n$ and $m$ stand for the number of cities (excluding the fictitious city) and the number of days in the campaign period, respectively.

The average objective values, CPU times and gaps in the first three tables have been calculated for the pool of 38 instances in which GUROBI was able to find a best feasible solution (BFS) within a time limit of 24 hours. The BFS constitutes the tightest lower bound on the true optimal objective value that can be obtained by GUROBI within the specified time limit. For the remaining 7 instances of `PE.I` and `PE.II`, we benchmarked MS-IVND with FDOR only. The results of these latter tests are presented in Table 5. Proven optimal objective values are indicated by an asterisk (*), while the objective value (Obj. Val.) of the best known solution (BKS) for a specific test instance is shown in **boldface**.

Table 2. Results for 18 `PE.I` instances solvable by GUROBI

| `PE.I` Instances | GUROBI | | | | FDOR | | | MS-IVND | | |
|---|---|---|---|---|---|---|---|---|---|---|
| | BFS | Opt.Gap (%) | $t^{\text{BFS}}$ (s) | CPU (s) | Obj. Val. | CPU (s) | FDOR Gap (%) | Obj. Val. | CPU (s) | MS-IVND Gap (%) |
| `6C-2D` | **7110*** | 0.0 | 0 | 0.1 | **7110** | 0.1 | 0.0 | **7110** | 109.6 | 0.0 |
| `6C-3D` | **8181*** | 0.0 | 0 | 0.1 | **8181** | 0.1 | 0.0 | **8181** | 138.0 | 0.0 |
| `7C-2D` | **9629*** | 0.0 | 0 | 0.2 | **9629** | 0.1 | 0.0 | **9629** | 124.4 | 0.0 |
| `7C-4D` | **11597*** | 0.0 | 0 | 0.4 | 11457 | 0.2 | 1.2 | **11597** | 135.2 | 0.0 |
| `9C-3D` | **10939*** | 0.0 | 0 | 0.5 | 10788 | 0.1 | 1.4 | **10939** | 150.7 | 0.0 |
| `9C-4D` | **11668*** | 0.0 | 1 | 1.3 | 11268 | 0.1 | 3.4 | **11668** | 183.0 | 0.0 |
| `12C-5D` | **14575*** | 0.0 | 6 | 6.0 | 12906 | 0.3 | 11.5 | 13682 | 207.6 | 6.1 |
| `15C-7D` | **17240*** | 0.0 | 462 | 551.3 | 16132 | 0.5 | 6.4 | 16491 | 209.2 | 4.3 |
| `15C-10D` | **18759*** | 0.0 | 19972 | 30458.5 | 17356 | 0.7 | 7.5 | 18065 | 22.8.0 | 3.7 |
| `21C-7D` | **19138*** | 0.0 | 2026 | 6705.3 | 17325 | 0.9 | 9.5 | 17709 | 263.8 | 7.5 |
| `21C-10D` | **21904** | 6.8 | 11582 | 86400.0 | 20673 | 1.2 | 5.6 | 20934 | 230.3 | 4.4 |
| `30C-7D` | **29427*** | 0.0 | 20665 | 20670.3 | 27474 | 1.7 | 6.6 | 28582 | 217.8 | 2.9 |
| `30C-10D` | **35013** | 5.9 | 30040 | 86400.0 | 32213 | 2.2 | 8.0 | 33109 | 286.9 | 5.4 |
| `40C-7D` | **30086** | 4.0 | 59757 | 86400.0 | 28821 | 3.7 | 4.2 | 29187 | 213.0 | 3.0 |
| `40C-10D` | **36409** | 12.6 | 62342 | 86400.0 | 34672 | 4.9 | 4.8 | 36005 | 392.0 | 1.1 |
| `51C-7D` | **41087** | 9.9 | 85597 | 86400.0 | 36942 | 8.4 | 10.1 | 37449 | 423.1 | 8.9 |
| `51C-10D` | **45667** | 22.3 | 77316 | 86400.0 | 43212 | 11.3 | 5.4 | 44047 | 484.7 | 3.6 |
| `51C-30D` | 47279 | 186.7 | 61189 | 86400.0 | 59890 | 14.5 | –26.7 | **63009** | 829.7 | –33.3 |
| *Average* | 23094.9 | | 23941.9 | 36844.1 | 22558.3 | 2.8 | 3.27 | **23188.5** | 270.5 | 0.98 |

Table 3. Results for 17 `PE.II` instances solvable by GUROBI

| `PE.II` Instances | GUROBI | | | | FDOR | | | MS-IVND | | |
|---|---|---|---|---|---|---|---|---|---|---|
| | BFS | Opt.Gap (%) | $t^{BFS}$ (s) | CPU (s) | Obj. Val. | CPU (s) | FDOR Gap (%) | Obj. Val. | CPU (s) | MS-IVND Gap (%) |
| `20C-5D` | **25118*** | 0.0 | 44 | 239.2 | 24196 | 0.6 | 3.7 | **25118** | 306.3 | 0.0 |
| `20C-7D` | **27523*** | 0.0 | 454 | 1995.9 | 25419 | 0.6 | 7.6 | **27523** | 440.7 | 0.0 |
| `30C-5D` | **16635*** | 0.0 | 161 | 709.9 | 16052 | 1.5 | 3.5 | **16635** | 212.0 | 0.0 |
| `30C-7D` | **18855*** | 0.0 | 13163 | 28216.8 | 17997 | 1.8 | 4.6 | **18855** | 228.6 | 0.0 |
| `30C-10D` | **21251** | 5.9 | 17722 | 86400.0 | 19577 | 2.0 | 7.9 | 20229 | 219.8 | 4.8 |
| `40C-7D` | 32811 | 20.1 | 76679 | 86400.0 | 31748 | 3.0 | 3.2 | **33023** | 266.8 | –0.6 |
| `40C-10D` | **37851** | 3.8 | 51297 | 86400.0 | 34267 | 3.6 | 9.5 | 36427 | 351.7 | 3.7 |
| `50C-7D` | 32829 | 1.6 | 4945 | 86400.0 | **33101** | 6.9 | –0.8 | 32881 | 321.3 | –0.2 |
| `50C-10D` | **38098** | 11.8 | 45006 | 86400.0 | 37389 | 8.5 | 1.9 | 37344 | 325.0 | 2.0 |
| `50C-15D` | **44098** | 35.6 | 70662 | 86400.0 | 41687 | 11.0 | 5.5 | 42293 | 575.4 | 4.1 |
| `60C-7D` | **40480** | 2.5 | 36955 | 86400.0 | 38105 | 13.8 | 5.9 | 40237 | 217.7 | 0.6 |
| `60C-10D` | **48270** | 7.0 | 82709 | 86400.0 | 45446 | 18.8 | 5.8 | 48066 | 385.9 | 0.4 |
| `60C-20D` | 50559 | 80.1 | 64244 | 86400.0 | 62869 | 22.8 | –24.3 | **64056** | 565.1 | –26.7 |
| `70C-10D` | 42474 | 13.9 | 82434 | 86400.0 | 40201 | 26.1 | 5.4 | **45159** | 237.9 | –6.3 |
| `70C-20D` | 43705 | 112.3 | 83589 | 86400.0 | 51055 | 34.9 | –16.8 | **57439** | 591.3 | –31.4 |
| `80C-10D` | 40808 | 22.2 | 52003 | 86400.0 | 38423 | 38.6 | 5.8 | **42559** | 255.3 | –4.3 |
| `80C-20D` | 50777 | 75.1 | 74448 | 86400.0 | 53270 | 41.9 | –4.9 | **55691** | 350.0 | –9.7 |
| *Average* | 36008.4 | | 44500.9 | 67903.6 | 35929.5 | 13.9 | 1.39 | **37855.0** | 344.2 | –3.74 |

Table 4. Results for the `LE` instances

| `LE` Instances | GUROBI | | | | FDOR | | | MS-IVND | | |
|---|---|---|---|---|---|---|---|---|---|---|
| | BFS | Opt.Gap (%) | $t^{BFS}$ (s) | CPU (s) | Obj. Val. | CPU (s) | FDOR Gap (%) | Obj. Val. | CPU (s) | MS-IVND Gap (%) |
| `39C-7D` | 22361 | 4.7 | 85841 | 86400 | 22164 | 11.7 | 0.9 | **23360** | 308.2 | –4.5 |
| `39C-10D` | 26774 | 18.5 | 63714 | 86400 | 27191 | 13.5 | −1.6 | **27241** | 276.2 | –1.7 |
| `39C-14D` | 30214 | 57.5 | 76094 | 86400 | 31757 | 15.9 | −5.1 | **32486** | 339.2 | –7.5 |
| *Average* | 26449.7 | | 75216.3 | 86400.0 | 27037.3 | 13.7 | –1.93 | **27576.3** | 307.9 | –4.57 |

Table 5. Results for 7 `PE.I` and `PE.II` instances unsolvable by GUROBI

| `PE.I` Instances | FDOR | | MS-IVND | | |
|---|---|---|---|---|---|
| | Obj. Val. | CPU (s) | Obj. Val. | CPU (s) | MS-IVND Gap (%) |
| `70C-15D` | 46818 | 16.6 | **47813** | 304.7 | –2.1 |
| `70C-40D` | 58408 | 22.2 | **66819** | 1244.3 | –14.4 |
| `93C-30D` | 68174 | 26.6 | **69444** | 1191.6 | –1.9 |
| `93C-40D` | 73574 | 27.1 | **75148** | 1798.9 | –2.1 |
| `PE.II` Instances | | | | | |
| `70C-30D` | 57065 | 36.2 | **61786** | 994.1 | –8.3 |
| `80C-30D` | 57285 | 50.0 | **61605** | 738.4 | –7.5 |
| `80C-40D` | 62576 | 48.7 | **65027** | 1338.4 | –3.9 |
| *Average* | 60557.1 | 32.5 | **63948.9** | 1087.2 | –5.74 |

When GUROBI attains a proven optimal solution on a given instance, the associated optimality gap between the final lower and upper bounds of the GUROBI solution with respect to the lower bound BFS drops to 0.00%. Optimality gaps reported by GUROBI are shown under the column header “Opt.Gap (%)”. CPU times are measured in seconds. The columns with the header $t^{\text{BFS}}$ in Table 2, Table 3 and Table 4 reveal the time it took GUROBI to attain its BFS. The percent gaps in the first three tables are given under the column header “MS-IVND Gap (%)”. They have been calculated with respect to the BFS using the formula $\text{Gap (\%)} = 100 \times \frac{\text{BFS} - \text{Obj.Val.}}{\text{BFS}}$. On the other hand, the percent gaps in Table 5 have been calculated for the incumbent objective values of MS-IVND with respect to those of FDOR since GUROBI cannot attain a feasible solution despite running for 24 hours on the instances of Table 5.

In terms of solution quality, MS-IVND outperforms GUROBI and FDOR on 11 and 40 instances, respectively. There is a tie between GUROBI and MS-IVND in 10 instances. Nonetheless, the performance of MS-IVND is timewise far superior to GUROBI. With reference to the BFS of GUROBI, it achieves a smaller average gap (0.97%, –3.74%, –4.57%) than FDOR (3.27%, 1.39%, –1.93%) in all three instance sets. The average Obj. Val. is improved by 0.41%, 5.13%, and 4.71% compared with GUROBI, and by 2.79%, 5.36%, and 2.43% compared with FDOR. The percent improvement over FDOR is 5.60% in the hardest 7 instances which are unsolvable by GUROBI. MS-IVND finds 7 optimal solutions more than FDOR. In addition, it achieves 17 new BKSs.

## 5.4 Starting GUROBI at the incumbent solution of MS-IVND

Another set of experiments with GUROBI is performed as follows. We selected all 30 MPTPP instances where GUROBI was previously unable to find an optimal solution. For each instance we retrieved the values of the decision variables of the implied MILP model from the incumbent MS-IVND solution, and passed them to GUROBI as the initial decision variables vector using the `MIPStart` attribute. Subsequently, we started GUROBI at this incumbent MS-IVND solution. The new BFS, Opt.Gap (%) and CPU time results are reported and compared with the original results in Table 6. Note that the initial decision variables vector of an optimization model in GUROBI is set to zero unless instructed otherwise.

Starting GUROBI at the incumbent solution of MS-IVND led to one more optimal solution; instance `PE.I 40C7D` is now solved to proven optimality. The final BFS of GUROBI improved in a total of 18 instances. In three instances it remained the same. Interestingly, in two instances (`PE.I 51C7D` and `PE.II 50C15D`), the BFS deteriorated which means that the best possible net benefit value decreased. It is hard to comment on this retrogression as we are not exactly knowledgeable about the internal algorithms of the black box commercial solvers. It can also be linked to the structure of the new initial solution. Moreover, in 10 instances GUROBI achieved no improvement over the respective MS-IVND solution which was fed as the initial solution.

GUROBI documentation (GUROBI Optimization, 2020) says that starting the solver at an initial feasible nonzero solution will not always lead to a better solution at the end. All in all, this approach improved the solution quality in 30 test instances, but not dramatically. Average relative BFS improvement is about 6.58% in those instances where GUROBI was able to return a feasible solution previously.

Table 6. GUROBI solutions before and after starting with to the incumbent MS-IVND solution.

| Instances | Original GUROBI Solution | | | MS-IVND solution | GUROBI solution using the MS-IVND solution as the initial solution | | |
|---|---|---|---|---|---|---|---|
| **PE.I** | BFS | Opt. Gap (%) | CPU (s) | Obj. Val. | BFS | Opt. Gap (%) | CPU (s) |
| 21C-10D | 21904 | 6.8 | 86400 | 20934 | 21904 | 6.8 | 86400 |
| 30C-10D | 35013 | 5.9 | 86400 | 33109 | 35034 | 5.9 | 86400 |
| 40C-7D | 30086 | 4.0 | 86400 | 29187 | **30122** | 0.0 | 64501 |
| 40C-10D | 36409 | 12.6 | 86400 | 36005 | 36560 | 12.4 | 86400 |
| 51C-7D | 41087 | 9.9 | 86400 | 37449 | 38838 | 10.8 | 86400 |
| 51C-10D | 45667 | 22.3 | 86400 | 44047 | 46319 | 19.6 | 86400 |
| 51C-30D | 47279 | 186.7 | 86400 | 63009 | 63009 | 94.1 | 86400 |
| 70C-15D | – | – | 86400 | 47813 | 47905 | 67.3 | 86400 |
| 70C-40D | – | – | 86400 | 66819 | 66819 | 57.9 | 86400 |
| 93C-30D | – | – | 86400 | 69444 | 69822 | 96.1 | 86400 |
| 93C-40D | – | – | 86400 | 75148 | 75164 | 52.7 | 86400 |
| **PE.II** | | | | | | | |
| 30C-10D | 21251 | 5.9 | 86400 | 20229 | 21606 | 3.7 | 86400 |
| 40C-7D | 32811 | 2.1 | 86400 | 33023 | 33026 | 1.5 | 86400 |
| 40C-10D | 37851 | 3.8 | 86400 | 36427 | 37851 | 3.8 | 86400 |
| 50C-7D | 32829 | 1.6 | 86400 | 32881 | 33319 | 0.1 | 86400 |
| 50C-10D | 38098 | 11.8 | 86400 | 37344 | 38235 | 11.5 | 86400 |
| 50C-15D | 44098 | 35.6 | 86400 | 42293 | 42293 | 32.6 | 86400 |
| 60C-7D | 40480 | 2.5 | 86400 | 40237 | 40480 | 2.5 | 86400 |
| 60C-10D | 48270 | 7.0 | 86400 | 48066 | 48515 | 6.2 | 86400 |
| 60C-20D | 50559 | 80.1 | 86400 | 64056 | 64056 | 42.8 | 86400 |
| 70C-10D | 42474 | 13.9 | 86400 | 45159 | 45159 | 10.1 | 86400 |
| 70C-20D | 43705 | 112.3 | 86400 | 57439 | 57458 | 78.5 | 86400 |
| 70C-30D | – | – | 86400 | 61786 | 61965 | 59.6 | 86400 |
| 80C-10D | 40808 | 22.2 | 86400 | 42559 | 42559 | 17.9 | 86400 |
| 80C-20D | 50777 | 75.1 | 86400 | 55691 | 58053 | 88.0 | 86400 |
| 80C-30D | – | – | 86400 | 61605 | 61630 | 58.6 | 86400 |
| 80C-40D | – | – | 86400 | 65027 | 65027 | 64.9 | 86400 |
| **LE** | | | | | | | |
| 39C-7D | 22361 | 4.7 | 86400 | 23360 | 23360 | 0.4 | 86400 |
| 39C-10D | 26774 | 18.5 | 86400 | 27241 | 27241 | 16.0 | 86400 |
| 39C-14D | 30214 | 57.5 | 86400 | 32486 | 32486 | 47.5 | 86400 |
| *Average* | 37426.3 | 30.56 | 86400 | 39227.4 | 39890.6 | 22.29 | 85448 |

### 5.5 Results on a Case Study

We next compare the real-life campaign plan of the party with the MILP model solution and the best heuristic solution found by MS-IVND. This comparison better emphasises the need to solve the MPTPP. To this end, we retrieved the opposition party's realised meetings prior to the general election in June 2015. In the light of these meetings, we created our large-size instance with 70 cities and a campaign period of 40 days. In order to make a fair comparison, we removed from the MPTPP model the constraints associated with two of our assumptions.

(i) The first constraint was forcing the politician to hold at least one meeting every day. However, in the actual meeting schedule of the party there were two meeting-free days.

(ii) The second constraint was forcing the politician to end the campaign at the campaign centre. We also lifted this constraint since the actual campaign of the party back in June 2015 had not been completed in Ankara. Table 7 illustrates the solutions of GUROBI and MS-IVND alongside the party's actual plan on the `70C-40D` instance.

Table 7 Comparison of the GUROBI and MS-IVND solutions with the actual schedule of the party

| | Obj. Val. | Opt. Gap(%) | # of Meetings | CPU time |
|---|---|---|---|---|
| GUROBI | LB = 46,640 (BFS)<br>UB = 117,427 | 60.3 | 75 | 259,200 s<br>(3 days) |
| Party's Plan | 24,534 | — | 77 | — |
| MS-IVND | 64,830 | — | 98 | 1012 s |

Table 7 shows that GUROBI is not able to solve the new MPTPP model to optimality despite a time allowance of 3 days. It can only bracket the true optimal objective value between a lower bound (LB) of 46,640 and an upper bound (UB) of 117,427 where the former is the objective value of the BFS reached by GUROBI. However, the BFS bears a net benefit that is about 90% greater than the net benefit accrued by the end of the actual campaign plan of the party. In the actual plan there are three meetings in Istanbul, Ankara, and Mersin each; two meetings in Izmir; and one meeting in the remaining cities. However, the BFS of GUROBI prescribes three meetings in İstanbul, Ankara, İzmir, and Mersin each; two meetings in the majority of midsize cities such as Adana, Balıkesir, Bursa, Çanakkale, Hatay, Konya, Zonguldak, Uşak, etc.; and one meeting in the remaining cities. MS-IVND, on the other hand, is able to find a much better solution with 98 meetings. The results in Table 7 underline the massive advantage of solving the MPTPP for the maximisation of the net benefit obtained from an election campaign that involves a relatively large number of cities and spans an extended period.

# 6 Scenario analysis and managerial insights

In this section we conduct an extensive scenario analysis to gain managerial insights into the MPTPP. We consider two levels of scenarios; the first level covers extreme scenarios, and the second one investigates the effect of different objective functions. Scenario descriptions are presented in the next two subsections followed by the last one which discusses their respective results.

## 6.1 Scenario analysis level 1: Extreme scenarios

We consider the following four scenarios with $\phi = 14$ hours as the maximum tour duration (MTD).

*Scenario* 1: Base Scenario, described in section 3.

*Scenario* 2: Base Scenario with the additional restriction of at most one meeting in each candidate city throughout the campaign.

*Scenario* 3: Base Scenario where the objective function involves only collected rewards and no travelling costs. In addition, the politician needs not return to the capital city periodically.

*Scenario* 4: Base Scenario with only closed daily tours originating and terminating at the capital city every day.

### 6.1.1 Description of the scenarios

*Scenario 1: Model Full-MILP (Base Scenario)*

The politician's campaign starts in Ankara. The politician cannot be away from Ankara for more than $\kappa = 5$ days in a row. There is no restriction as to where to terminate the tour (sleep) at the end of a given day and start the tour of the next day (wake up). Thus, the tour on a given day $t$ can be either an open or a closed tour. There is an MTD constraint in place which prohibits daily tours in excess of 14 hours. Candidate cities are divided into three groups, namely big cities (İstanbul, Ankara and İzmir), midsize cities and small cities where the number of meetings hosted is limited to three, two and one, respectively. The number of meetings held each day is limited to four.

*Scenario 2: Model Full-1Meet*

This scenario is derived from Scenario 1 by revoking the option of multiple meetings in big and midsize cities. Since each city can host at most one meeting during the campaign, the binary decision variables $R_{its}$ are excluded from the MILP model of the Base Scenario which, in turn, simplifies the model of the MPTPP drastically. We then modify our MS-IVND according to the no-repeated-meeting assumption. While exploring the neighbourhoods, a list of all meetings is kept, and the algorithm performs another feasibility check to avoid repeated meetings. The net benefit definition comprising the objective function is simplified as shown in (5). It is also worth noting that the optimal solutions (thus the optimal objective values) of Scenario 1 and Scenario 2 may be identical.

$$NET\ BENEFIT(Full\text{-}1Meet) = \sum_{i\in N}\sum_{t\in T}\pi_i \frac{m-t+1}{m} Z_{it} - \bar{K}\sum_{i\in N}\sum_{j\in N}\sum_{t\in T} c_{ij} X_{ijt} \qquad (5)$$

*Scenario* 3*: Model Rew-Only*

The third scenario is derived from *Scenario 1* using the following two modifications:

(i) The necessity to periodically return to the capital Ankara at least once every $\kappa$ days is lifted. The politician has full freedom to hop from one city to another. Also, the politician can stay overnight in any city. Yet the campaign is still going to start in Ankara on day 1.

(ii) Travelling costs are discarded from the objective function. This fundamental change motivates the politician to roam between all candidate cities without considering the cost of travelling.

These two modifications make our MS-IVND relatively much faster, because the computationally expensive feasibility check of returning to the campaign centre is now redundant.

*Scenario 4: Model Alt-1Depot*

The fourth scenario is derived from Scenario 1 by a fundamental paradigm shift in which the politician wakes up in the capital city Ankara every morning and returns there to sleep by the end of every day. This implies that each daily trip is going to be a closed tour with Ankara being the depot of the trip. This reduces the problem to a multi-period selective TSP with a single depot. In our MS-IVND implementation, we do not permit the algorithm to apply the neighbourhood moves on depot nodes. Note that *Alt-1Depot* is a much more restrictive model than *Full-MILP* since it does not allow open tours. Clearly, the optimal objective value of *Alt-1Depot* is a valid lower bound on that of *Full-MILP*. *Alt-1Depot* is also computationally much more tractable due to the absence of the novel binary variables representing exclusive departures from and exclusive entrances to candidate cities. The reader is referred to Shahmanzari et al. (2020) for an in-depth discussion of those variables of the model *Full-MILP*.

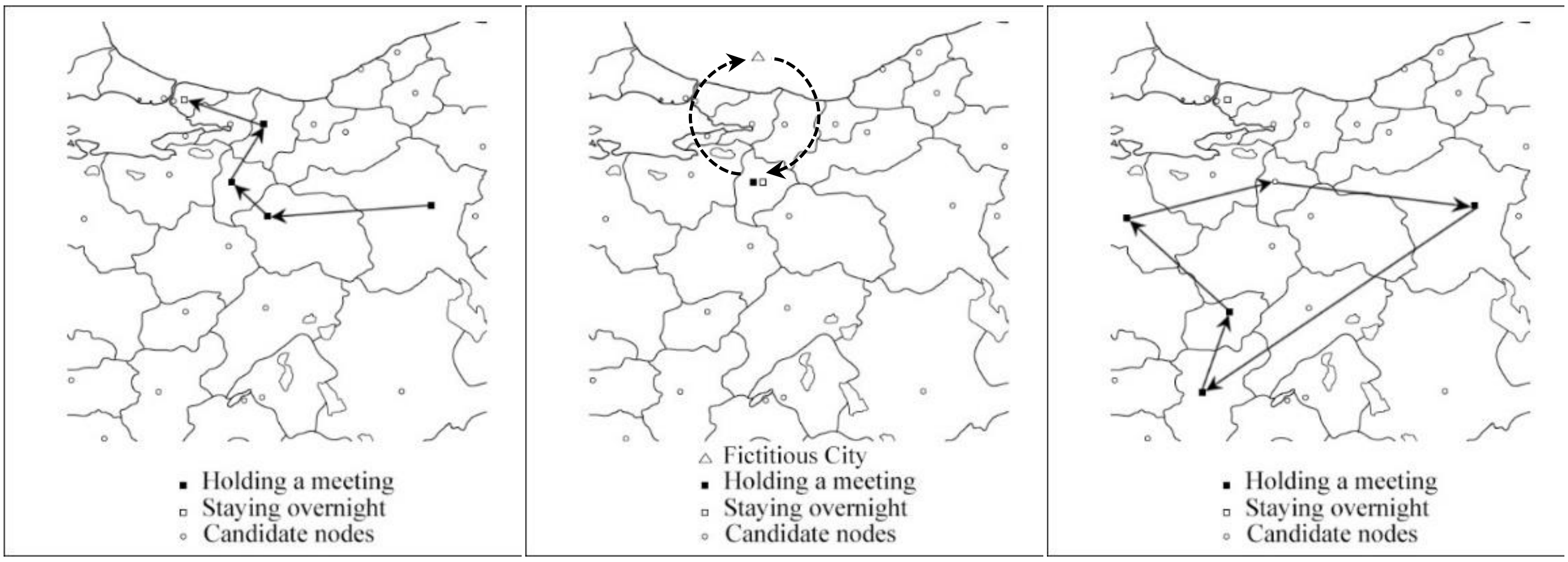

Fig. 6 Different tours in the first three scenarios.

### 6.2 Scenario analysis level 2: Alternative reward function

Thus far we assumed that earlier meetings produce higher rewards. In an alternative scenario, we reverse the direction of the reward function such that the reward increases as we approach the election day, i.e. the end of the campaign period. In this case, the objective function is defined in Eq. (6) where linearity of the objective function is also preserved.

$$\max.\ (NET\ BENEFIT) = \sum_{i\in\mathbf{N}}\sum_{t\in\mathbf{T}}\frac{t+m}{m}\pi_i\, FM_{it} + \sum_{i\in\mathbf{N}}\sum_{t\in\mathbf{T}}\sum_{1\le s<t}\frac{s(t+m)}{Km^2}\pi_i\, R_{its} - \sum_{i\in\mathbf{N}}\sum_{j\in\mathbf{N}}\sum_{t\in\mathbf{T}}c_{ijt}X_{ijt} \qquad (6)$$

### 6.3 Computational results of the scenario analysis

#### 6.3.1 Scenario analysis level 1

The results of the four different scenarios of level 1 are presented in Table 8. The naming convention in the leftmost column of the table sheds light on the sizes of the 10 test instances. We choose only 10 instances to illustrate the characteristics of each scenario where we ensure that there exist small-sized, medium-sized, and large-sized cities in every instance. **Boldface** figures in the BFS column of each scenario point to proven optimality achieved by the commercial solver GUROBI. Possible tours in the first three scenarios are illustrated in Fig. 6. Table 8 consists of four segments where in each segment the first three columns correspond to GUROBI. These columns show the optimal or best feasible objective value, the number of meetings held during the planning horizon $(\alpha')$, and the final gap attained by GUROBI. Within each segment, the second group of three columns correspond to the MS-IVND results, and show the objective value and the final gap of MS-IVND.

The CPU time limit in all GUROBI runs was 24 hours. The average optimality gap is 6.19% which implies MPTPP is a large-scale optimisation problem even for small-size instances. In *Full-1Meet* (Scenario 2), the computational complexity of the problem is greatly reduced due to the removal of the binary variables of repeated meetings from the model and due to the simplified net benefit definition shown in Eq (5). However, the average optimality gap of GUROBI is still high (7.20%) in Scenario 2.

Objective values obtained in *Rew-Only* (Scenario 3) are not comparable with the ones obtained in the other scenarios since the travelling costs in the definition of net benefit are ignored. However, except in three instances, namely `40C-7D`, `51C-7D` and `51C-10D`, the count of meetings realised in this scenario is either higher than or equal to the count of meetings in other scenarios. This can be ascribed to having a larger feasible solution space which occurs because of lifting the necessity to visit the campaign centre every $\kappa$ days. Similarly, the politician in this scenario has more freedom to travel to remote cities that would not be visited in the Base Scenario due to the net benefit being negative after the deduction of travelling expenses. The average gap in Scenario 3 is 6.07%.

Table 8. Results of the four scenarios

| | Scenario 1: *Full-MILP* | | | | | | Scenario 2: *Full-1Meet* | | | | | |
|---|---|---|---|---|---|---|---|---|---|---|---|---|
| | GUROBI | | | MS-IVND | | | GUROBI | | | MS-IVND | | |
| Instance | BFS | $\alpha'$ | Opt.Gap (%) | Obj. Val. | $\alpha'$ | MS-IVND Gap (%) | BFS | $\alpha'$ | Opt.Gap (%) | Obj. Val. | $\alpha'$ | MS-IVND Gap (%) |
| 15C-7D | **17240** | 14 | 0.0 | 16491 | 15 | 4.3 | **16000** | 13 | 0.0 | **16000** | 13 | 0.0 |
| 15C-10D | **18759** | 17 | 0.0 | 18065 | 17 | 3.7 | **16299** | 15 | 0.0 | **16299** | 15 | 0.0 |
| 21C-7D | **19138** | 16 | 0.0 | 17709 | 17 | 7.5 | **18117** | 16 | 0.0 | **18117** | 16 | 0.0 |
| 21C-10D | 21904 | 21 | 6.9 | 20934 | 19 | 4.4 | **20850** | 20 | 0.0 | **20850** | 20 | 0.0 |
| 30C-7D | **29427** | 18 | 0.0 | 28582 | 18 | 2.9 | 27576 | 17 | 3.3 | 27912 | 23 | –1.2 |
| 30C-10D | 35013 | 24 | 6.0 | 33109 | 25 | 5.4 | 32210 | 24 | 9.0 | 33964 | 26 | –5.4 |
| 40C-7D | 30086 | 20 | 4.1 | 29187 | 18 | 3.0 | 27023 | 18 | 16.4 | 28550 | 19 | –5.6 |
| 40C-10D | 36409 | 25 | 12.6 | 36005 | 26 | 1.6 | 32210 | 24 | 9.0 | 35021 | 28 | –8.7 |
| 51C-7D | 41087 | 31 | 9.9 | 37449 | 30 | 8.9 | 33366 | 23 | 17.6 | 37042 | 25 | –11.0 |
| 51C-10D | 45667 | 36 | 22.4 | 44047 | 35 | 3.5 | 35314 | 26 | 16.7 | 40117 | 30 | –13.6 |

| | Scenario 3: *Rew- Only* | | | | | | Scenario 4: *Alt-1Depot* | | | | | |
|---|---|---|---|---|---|---|---|---|---|---|---|---|
| | GUROBI | | | MS-IVND | | | GUROBI | | | MS-IVND | | |
| Instance | BFS | $\alpha'$ | Opt.Gap (%) | Obj. Val. | $\alpha'$ | MS-IVND Gap (%) | BFS | $\alpha'$ | Opt.Gap (%) | Obj. Val. | $\alpha'$ | MS-IVND Gap (%) |
| 15C-7D | **22561** | 16 | 0.0 | 21730 | 16 | 3.7 | **11539** | 10 | 0.0 | **11539** | 10 | 0.0 |
| 15C-10D | **26061** | 20 | 0.0 | 25238 | 17 | 3.2 | 11170 | 13 | 15.4 | 11834 | 13 | –5.9 |
| 21C-7D | **23932** | 17 | 0.0 | 21750 | 16 | 9.1 | **13779** | 12 | 0.0 | **13779** | 12 | 0.0 |
| 21C-10D | 28927 | 23 | 3.8 | 26947 | 21 | 6.8 | 14498 | 15 | 14.6 | 15013 | 15 | –3.6 |
| 30C-7D | **33638** | 18 | 0.0 | 33115 | 19 | 1.6 | **20774** | 13 | 0.0 | **20774** | 13 | 0.0 |
| 30C-10D | 38148 | 24 | 10.0 | 36991 | 24 | 3.0 | 24520 | 17 | 15.9 | 25901 | 18 | –5.6 |
| 40C-7D | 32586 | 19 | 9.5 | 31435 | 19 | 3.5 | 20893 | 13 | 7.1 | 21375 | 13 | –2.3 |
| 40C-10D | 39876 | 25 | 13.7 | 39226 | 27 | 1.6 | 25324 | 18 | 19.73 | 26230 | 19 | –3.6 |
| 51C-7D | 33123 | 19 | 11.0 | 30774 | 16 | 7.1 | **31154** | 19 | 0.0 | **31154** | 19 | 0.0 |
| 51C-10D | 41325 | 28 | 12.7 | 38827 | 25 | 6.0 | 37373 | 26 | 5.1 | 37824 | 29 | –1.2 |

The results of *Alt-1Depot* (Scenario 4) are also interesting, as this scenario bears the most similar conditions to the current campaign policy of the party. In comparison to *Full-MILP*, the commercial solver GUROBI was able to attain optimality in one more instance (51C-7D). Table 8 reports an overall lesser number of meetings in Scenario 4. It is apparent that the requirement to return to the capital city Ankara at the end of every day prevents some of the meetings which were realised in the *Full-MILP* model of the Base Scenario. The average gap reported by GUROBI in this scenario is 7.78%.

### 6.3.2 Scenario analysis level 2

The importance of holding meetings in early days or in the last days of the campaign period should be decided by party executives indeed. Despite this fact, we present in Table 9 the comparison between the original and the alternative reward functions on 14 small-size test instances where we solved the *Full-MILP* model of the Base Scenario. All instances are solved to proven optimality under each net benefit function. The column with the header “CPU (s)” indicates the solution times in seconds reported by GUROBI.

According to Table 9 the solution times obtained with the original reward function are better in 13 out of 14 instances. The objective values are unfortunately not comparable due to different rewards being assigned to each city in different days. The decision to choose which reward function is obviously up to the politician. Based on the results in Table 9, we can say that the use of the alternative reward function in (6) as the objective function of the MPTPP increases the solution times of our algorithm MS-IVND.

Table 9 Comparison of new reward function and original reward function

| | *Full-MILP* with original reward function | | | *Full-MILP* with alternative reward function (6) | | |
|---|---|---|---|---|---|---|
| Instances | Obj. Value | Gap (%) | CPU (s) | Obj. Value | Gap (%) | CPU (s) |
| 6C-2D | 7110.8 | 0.0 | 0.1 | 17441.2 | 0.0 | 0.1 |
| 6C-3D | 8181.3 | 0.0 | 0.1 | 22272.2 | 0.0 | 0.2 |
| 7C-2D | 9629.5 | 0.0 | 0.1 | 22172.4 | 0.0 | 0.2 |
| 7C-3D | 10939.0 | 0.0 | 0.2 | 26778.4 | 0.0 | 0.4 |
| 7C-4D | 11597.4 | 0.0 | 0.4 | 30339.9 | 0.0 | 0.7 |
| 9C-2D | 9695.0 | 0.0 | 0.3 | 22172.4 | 0.0 | 0.1 |
| 9C-3D | 10939.0 | 0.0 | 0.5 | 28572.9 | 0.0 | 0.9 |
| 9C-4D | 11668.4 | 0.0 | 1.3 | 32149.5 | 0.0 | 1.8 |
| 12C-3D | 12620.0 | 0.0 | 1.1 | 31726.9 | 0.0 | 1.3 |
| 12C-4D | 13584.8 | 0.0 | 2.3 | 37076.3 | 0.0 | 3.6 |
| 12C-5D | 14575.6 | 0.0 | 6.0 | 40382.8 | 0.0 | 6.4 |
| 15C-3D | 12620.0 | 0.0 | 1.6 | 32750.8 | 0.0 | 2.9 |
| 15C-4D | 14210.3 | 0.0 | 4.3 | 39496.6 | 0.0 | 6.1 |
| 15C-5D | 15446.9 | 0.0 | 14.4 | 43533.7 | 0.0 | 159.8 |

# 7 Conclusions

In this paper we investigate a logistical planning problem arising in election campaigns which is known as the Multi-Period Travelling Politician Problem (MPTPP). It involves a politician who wants to maximise the net benefit of his/her campaign over a fixed period of days. Time-dependent and meeting frequency-dependent rewards are earned by holding a meeting in a city visited on a daily tour. The objective function represents the net benefit defined as the total collected reward minus the total travelling cost. Several real-life aspects such as the necessity to return to the campaign centre frequently, the maximum tour duration, and a time- as well as recency-dependent reward function are also incorporated into the model.

We propose a hybrid multi-start metaheuristic which integrates ILS and VND. It also leverages perturbation and local search schemes where three characteristic moves are built into the former and four in the latter. Our approach restarts with a new initial point whenever VND reaches a local optimum. Computational results on 45 test instances published in the literature reveal that the proposed metaheuristic, which we call MS-IVND, is an efficient algorithm for solving the MPTPP. It outperforms a recently published two-phase matheuristic called FDOR and the commercial MILP solver GUROBI in terms of solution quality and speed, respectively, which serve as two basic performance criteria. MS-IVND produces 7 optimal solutions and 17 new best known solutions. The superiority of the new algorithm MS-IVND over the MILP solution approach with GUROBI lies in its solution speed. On the other hand, it has a significant solution quality advantage over the matheuristic FDOR. We conclude that MS-IVND can help to reap higher net benefits from a multi-period election campaign.

In order to gain managerial insights into this problem, we also carry out a scenario analysis using the Base Scenario model referred to as *Full-MILP*. We conduct extensive experiments in which we demonstrate the power as well as the flexibility of MS-IVND which yields favourable results within relatively short computational times. Our experimental results provide useful insights into the planning of an election campaign in a real-world setting. This study will hopefully inspire other researchers to explore new avenues in the research of election logistics, touristic trip planning, and marketing campaigns.

## Acknowledgement

The authors are thankful to three anonymous referees and to the Associate Editor for their comments and suggestions that filled in the margins for improvement in the original draft. The fast and scrupulous review process of the *Journal of the Operational Research Society* during the ongoing hardship of the COVID-19 pandemic is deeply appreciated.

# Appendix A

We provide the description of the MILP formulation of the MPTPP in this appendix. The formulation is originally due to Aksen and Shahmanzari (2018). The reader is also referred to Shahmanzari et al. (2020) for an in-depth discussion of the model constituents, another alternative formulation for the satisfaction of the maximum tour duration, and also for valid inequalities which help tighten the original formulation.

*Index Sets:*

| | |
|---|---|
| $\mathbf{N} = \{0,...,n\}$ | Set $\mathbf{V}$ joined by city '0' which denotes a fictitious city with all associated costs, rewards and meeting duration being zero. |
| $\mathbf{V} = \mathbf{N} \setminus \{0\}$ | The set of cities to be considered for collecting rewards throughout the election campaign where city $i = 1$ denotes the campaign centre. |
| $\mathbf{T} = \{1,...,m\}$ | The set of $m$ days comprising the campaign duration. |

*Parameters:*

| | |
|---|---|
| $c_{ij}$ | Travelling cost from city $i$ to $j$ where $c_{ii} = 0$. |
| $d_{ij}$ | Travelling time from city $i$ to city $j$ where $d_{ii} = 0$. |
| $\pi_i$ | The base reward of city $i$. |
| $\sigma_i$ | The meeting duration in city $i$. |

$\alpha$ Maximum number of meetings allowed each day.

$T_{\max}$ Maximum tour duration (in hours) in each daily tour.

$\kappa$ Maximum number of consecutive days during which the politician is allowed to be away from the campaign centre.

$K$ The base reward depreciation coefficient applied in successive meetings held in the same city.

$\bar{K}$ Normalisation coefficient multiplied with the collected rewards to make travelling costs and daily rewards compatible.

*Decision Variables:*

$X_{ijt}$ Binary variable indicating if arc $(i, j)$ is traversed on day $t$ $(i, j \in \mathbf{N},\ t \in \mathbf{T})$ with $X_{iit} = 0.$

$L_{it}$ Binary variable indicating if the politician does not enter, but only leaves city $i$ in day $t$.

$E_{it}$ Binary variable indicating if the politician does not leave, but only enters city $i$ in day $t$.

$S_{it}$ Binary variable indicating if the politician stays overnight (sleeps) in city $i$ by the end of day $t$.

$Z_{it}$ Binary variable indicating if the politician holds a meeting in city $i$ on day $t$

$FM_{it}$ Binary variable indicating if the first meeting in city $i$ is held on day $t$.

$R_{its}$ Binary variable indicating if city $i$ accommodates two consecutive meetings on day $t$ and day

$R_{its}$ $(t-s)$ with no other activity in between. Since $1 \le s < t$, we have $R_{its} = 0$ for $t \le s \le m$.

$U_{it}$ A continuous nonnegative variable used in the Modified Miller-Tucker-Zemlin subtour elimination constraints. It is used to determine the order of visit for city $i$ on day $t$.

The MPTPP can be formulated as the following mixed-integer linear program:

$$\begin{aligned} \max.\ NET\ BENEFIT = & \sum_{i\in\mathbf{N}} \sum_{t\in\mathbf{T}} \pi_i \frac{m-t+1}{m} FM_{it} + \sum_{i\in\mathbf{N}} \sum_{t\in\mathbf{T}} \sum_{1\le s<t} \pi_i \frac{m-t+1}{m} \times \frac{s}{Km} R_{its} \\ & - \bar{K} \sum_{i\in\mathbf{N}} \sum_{j\in\mathbf{N}} \sum_{t\in\mathbf{T}} c_{ij} X_{ijt} \end{aligned} \tag{7}$$

Subject to:

$$\sum_{j\in\mathbf{N}} X_{ijt} \le 1 \qquad i \in \mathbf{N},\ t \in \mathbf{T} \tag{8}$$

$$\sum_{j\in\mathbf{N}} X_{jit} \le 1 \qquad i \in \mathbf{N},\ t \in \mathbf{T} \tag{9}$$

$$\sum_{i\in\mathbf{V}} Z_{it} \le \alpha \qquad t \in \mathbf{T} \tag{10}$$

$$\sum_{i\in\mathbf{V}} Z_{it} \ge 1 \qquad t \in \mathbf{T} \tag{11}$$

$$\sum_{i\in\mathbf{V}} \sigma_i Z_{it} + \sum_{i\in\mathbf{N}} \sum_{j\in\mathbf{N}} d_{ij} X_{ijt} \le T_{\max} \qquad t \in \mathbf{T} \tag{12}$$

$$FM_{i1} = Z_{i1} \qquad i \in \mathbf{V} \tag{13}$$

$$FM_{it} \le Z_{it} \qquad i \in \mathbf{V},\ t \in \mathbf{T} \setminus \{1\} \tag{14}$$

$$FM_{it} \le 1 - Z_{iu} \qquad i \in \mathbf{V},\ t \in \mathbf{T} \setminus \{1\},\ 1 \le u < t \tag{15}$$

$$\sum_{j\in\mathbf{N}} X_{ijt} - \sum_{j\in\mathbf{N}} X_{jit} = L_{it} - E_{it} \qquad i\in\mathbf{N},\ t\in\mathbf{T} \qquad (16)$$

$$L_{it} + E_{it} \le 1 \qquad i\in\mathbf{N},\ t\in\mathbf{T} \qquad (17)$$

$$\sum_{i\in\mathbf{N}} (L_{it} + E_{it}) \le 2 \qquad t\in\mathbf{T} \qquad (18)$$

$$S_{i(t-1)} \le S_{it} + \sum_{j\in\mathbf{N}} \frac{L_{jt}+E_{jt}}{2} \qquad i\in\mathbf{N},\ t\in\mathbf{T}\setminus\{1\} \qquad (19)$$

$$\sum_{j\in\mathbf{N}} \frac{L_{jt}+E_{jt}}{2} + S_{i(t-1)} \ge S_{it} \qquad i\in\mathbf{N},\ t\in\mathbf{T}\setminus\{1\} \qquad (20)$$

$$S_{i(t-1)} \le L_{it} + S_{it} \qquad i\in\mathbf{V},\ t\in\mathbf{T}\setminus\{1\} \qquad (21)$$

$$S_{0t} = 0 \qquad t\in\mathbf{T} \qquad (22)$$

$$S_{it} \ge X_{i0t} \qquad i\in\mathbf{V},\ t\in\mathbf{T} \qquad (23)$$

$$S_{i(t-1)} \ge X_{i0t} \qquad i\in\mathbf{V},\ t\in\mathbf{T}\setminus\{1\} \qquad (24)$$

$$X_{i0t} = X_{0it} \qquad i\in\mathbf{V},\ t\in\mathbf{T} \qquad (25)$$

$$E_{it} \le S_{it} \qquad i\in\mathbf{V},\ t\in\mathbf{T} \qquad (26)$$

$$S_{it} \le \sum_{j\in\mathbf{N}} X_{ij(t+1)} \qquad i\in\mathbf{V},\ 1\le t<m \qquad (27)$$

$$\sum_{i\in\mathbf{V}} S_{it} = 1 \qquad t\in\mathbf{T} \qquad (28)$$

$$\sum_{k=t}^{t+\kappa} S_{1k} \ge 1 \qquad 1\le t\le m-\kappa \qquad (29)$$

$$Z_{it} \le \sum_{j\in\mathbf{N}} X_{ijt} + E_{it} \qquad i\in\mathbf{V},\ t\in\mathbf{T} \qquad (30)$$

$$Z_{it} \le \sum_{j\in\mathbf{N}} X_{jit} + L_{it} \qquad i\in\mathbf{V},\ t\in\mathbf{T} \qquad (31)$$

$$(\alpha+1)S_{j(t-1)} + (\alpha+1)(1-X_{ijt}) + U_{jt} \ge U_{it} + 1 \qquad i,j\ (i\ne j)\in\mathbf{N},\ t\in\mathbf{T}\setminus\{1\} \qquad (32)$$

$$U_{it} \le \alpha+1 \qquad i\in\mathbf{N},\ t\in\mathbf{T} \qquad (33)$$

$$U_{it} \le \sum_{j\in\mathbf{N}} \sum_{k\in\mathbf{N}} X_{jkt} + 1 \qquad i\in\mathbf{N},\ t\in\mathbf{T} \qquad (34)$$

$$U_{it} \ge S_{i(t-1)} \qquad i\in\mathbf{N},\ t\in\mathbf{T}\setminus\{1\} \qquad (35)$$

$$(\alpha+1)(1-S_{i(t-1)}) + 1 \ge U_{it} \qquad i\in\mathbf{N},\ t\in\mathbf{T}\setminus\{1\} \qquad (36)$$

$$U_{it} \ge \sum_{j\in\mathbf{N}} X_{ijt} \qquad i\in\mathbf{N},\ t\in\mathbf{T} \qquad (37)$$

$$U_{it} \ge S_{it} + \sum_{j\in\mathbf{N}} X_{ijt} \qquad i\in\mathbf{N},\ t\in\mathbf{T} \qquad (38)$$

$$U_{it} \le (\alpha+1)\sum_{j\in\mathbf{N}} X_{ijt} + (\alpha+1)\sum_{j\in\mathbf{N}} X_{jit} \qquad i\in\mathbf{N},\ t\in\mathbf{T} \qquad (39)$$

$$R_{its} \le Z_{it} \qquad i\in\mathbf{N},\ 2\le t\le m,\ 1\le s<t \qquad (40)$$

$$R_{its} \le Z_{i(t-s)} \qquad i\in\mathbf{N},\ 2\le t\le m,\ 1\le s<t \qquad (41)$$

$$\sum_{k=t-s+1}^{t-1} Z_{ik} \le s(1-R_{its}) \qquad i\in\mathbf{N},\ 3\le t<m,\ 2\le s<t \qquad (42)$$

$$R_{its}=0 \qquad i\in\mathbf{N},\ t\in\mathbf{T},\ t\le s\le m \qquad (43)$$

$$R_{ius}\le 1-FM_{it} \qquad i\in\mathbf{V},\ t\in\mathbf{T}\setminus\{1\},\ t<u<m,\ u-t<s<u \qquad (44)$$

$$R_{its}\ge Z_{i(t-s)}+Z_{it}-\sum_{k=t-s+1}^{t-1} Z_k -1 \qquad i\in\mathbf{V},\ 3\le t\le m,\ 2\le s<t \qquad (45)$$

$$X_{ijt},\ L_{it},\ E_{it},\ S_{it},\ Z_{it},\ FM_{it},\ R_{its}\in\{0,1\} \text{ and } U_{it}\ge 0 \qquad (46)$$

The objective function (7) maximises the total net benefit. Constraints (8)–(9) reflect the selective routing characteristic of the MPTPP. Constraints (10)–(11) determine a lower and an upper bound on the number of meetings held on a given day. Constraints (12) specify the maximum tour duration restriction on each daily tour of the politician. Constraints (13)–(15) couple the decision variables $FM_{it}$ and $Z_{it}$. Constraints (16)–(28) couple the binary decision variables $X_{ijt}$, $S_{it}$, $E_{it}$, and $L_{it}$. Constraints (29) ensure frequent visits to the campaign centre, at least once every $\kappa$ days. The constraint sets (30) and (31) guarantee that no reward can be collected from unvisited cities. Constraint (32)–(39) serve as subtour elimination constraints. Constraints (40)–(45) couple the decision variable $R_{its}$, $FM_{it}$ and $Z_{it}$. Finally, constraints (46) are nonnegativity and binary constraints imposed on the decision variables.

# Appendix B

Table B.1 summarises the GUROBI options chosen in implementation. `MIPGap` is computed as $|BFS-BPS|/|BFS|\times 100\%$ where $BFS$ and $BPS$ stand for the best feasible and best possible solutions, namely the tightest lower and upper bounds in a maximisation problem, respectively. The CPU time limit (`TimeLimit`) is set to 86,400 seconds (24 hours). The options `Threads` and `Concurrentmip` turn on the multithreading (concurrent optimisation) capabilities of GUROBI. When `Threads` is set to zero, the computing load is distributed onto all available fourteen cores (28 threads) of the Intel Xeon® E5-2690 v4 processor. On the other hand, when `Concurrentmip` is set to three, the solver divides available threads evenly between three independent MILP solve operations and performs them in parallel. Optimisation terminates when the first solve operation completes. In order to compare multithreading options, we tested the performance of GUROBI under different concurrent optimisation configurations. We observed that `Concurrentmip = 3` outperforms other configurations and finds the best feasible solution as well as the best possible solution achieving thereby the smallest optimality gap. `NumericFocus` controls the degree to which the code attempts to detect and manage numerical issues. It is set to 3 since the right-hand side values of the constraints are relatively large in our model. The reader is referred to GUROBI User's Manual (2018) for a more thorough explanation of these options.

Table B.1 List of GUROBI specific options applied to all runs.

| GUROBI specific options used in Python codes | | | |
|---|---|---|---|
| `MIPGap` | `= 0.000` | `Threads` | `= 0` |
| `TimeLimit` | `= 86400` | `Concurrentmip` | `= 3` |
| `IterationLimit` | `= 1.e9` | `NumericFocus` | `= 3` |
| `NodeLimit` | `= 5.0e8` | `DualReductions` | `= 0` |
| `Nodefilestart` | `= 6.5` | `InfUnbdInfo` | `= 1` |

# Appendix C

The candidate cities considered in the test instance `PE.I 40C-10D` and their population, base reward, meeting duration and criticality factor data are shown in Table C.1. The cities are sorted in descending order of their base rewards.

Table C.1. City characteristics in `PE.I 40C-10D`

| City | Population in 2015 | Base Reward | Meeting Duration | Criticality Factor |
|---|---|---|---|---|
| İstanbul | 14,657,434 | 2,370 | 2 | 5 |
| Ankara | 5,270,575 | 1,505 | 2 | 5 |
| İzmir | 4,168,415 | 1,295 | 2 | 5 |
| Bursa | 2,842,547 | 1,040 | 1.5 | 5 |
| Hatay | 1,533,507 | 1,000 | 1.5 | 5 |
| İskenderun | 247,220 | 1,000 | 1 | 5 |
| Antalya | 2,288,456 | 935 | 1.5 | 5 |
| Alanya | 134,396 | 935 | 1 | 5 |
| Adana | 2,183,167 | 736 | 1.5 | 4 |
| Kahramanmaraş | 1,096,610 | 680 | 1.5 | 4 |
| Gaziantep | 1,931,836 | 675 | 1.5 | 3 |
| Denizli | 993,442 | 660 | 1 | 4 |
| Aydın | 1,053,506 | 660 | 1.5 | 4 |
| Kocaeli | 1,780,055 | 645 | 1.5 | 3 |
| Gebze | 357,743 | 645 | 1 | 3 |
| Muğla | 908,877 | 640 | 1 | 4 |
| Çorlu | 273,362 | 640 | 1 | 4 |
| Mersin | 1,745,221 | 630 | 1.5 | 3 |
| Ordu | 728,949 | 580 | 1.5 | 4 |
| Manisa | 1,380,366 | 570 | 1 | 4 |
| Balıkesir | 1,186,688 | 525 | 1.5 | 3 |
| Kastamonu | 372,633 | 500 | 1 | 4 |
| Edirne | 402,537 | 500 | 1 | 4 |
| Kars | 292,660 | 480 | 1 | 4 |
| Eskişehir | 826,716 | 465 | 1 | 3 |
| Erzincan | 222,918 | 460 | 1 | 4 |
| Afyon | 709,015 | 435 | 1 | 3 |
| Adıyaman | 602,774 | 420 | 1 | 2 |
| Diyarbakır | 1,654,196 | 410 | 1.5 | 2 |
| Çanakkale | 513,341 | 405 | 1 | 3 |
| Isparta | 421,766 | 375 | 1 | 3 |
| Giresun | 426,686 | 375 | 1 | 3 |
| Kayseri | 1,341,056 | 370 | 1.5 | 2 |
| Konya | 2,130,544 | 362 | 1.5 | 2 |
| Amasya | 322,167 | 360 | 1 | 3 |
| Bolu | 291,095 | 360 | 1 | 3 |
| Niğde | 346,114 | 360 | 1 | 3 |
| Bartın | 190,708 | 330 | 1 | 3 |
| Malatya | 772,904 | 300 | 1 | 2 |
| Kırşehir | 225,562 | 230 | 1 | 2 |

The shaded areas in Fig. C.1 denote those cities of Turkey which are included in the 10-day-long campaign period. Full-MILP results on `PE.I 40C-10D` are provided in Fig. C.2 where the circled numbers on a city represent the number of meetings realised and the meeting days. Table C.2 reveals the daily tours where (M) indicates a meeting. The best gap obtained for this instance in 24 hours was 12.6%.

Table C.2. Daily tours of the instance **`PE.I 40C-10D`**

| Days | Route |
|---|---|
| 1 | Ankara (M) → Hatay (M) → İskenderun (M) |
| 2 | İskenderun → Adana (M) → Istanbul (M) |
| 3 | Istanbul → Kocaeli (M) → Bursa (M) → Balıkesir (M) |
| 4 | Balıkesir → Manisa (M) → İzmir (M) → Aydin (M) → Muğla |
| 5 | Muğla (M) → Denizli (M) → Antalya (M) → Isparta |
| 6 | Isparta (M) → Afyonkarahisar (M) → Eskişehir (M) → Ankara |
| 7 | Ankara (M) → Gebze (M) → Istanbul |
| 8 | Istanbul (M) → Gaziantep (M) → Kahramanmaraş (M) |
| 9 | Kahramanmaraş → Hatay (M) → Adana (M) → Mersin |
| 10 | Mersin (M) |

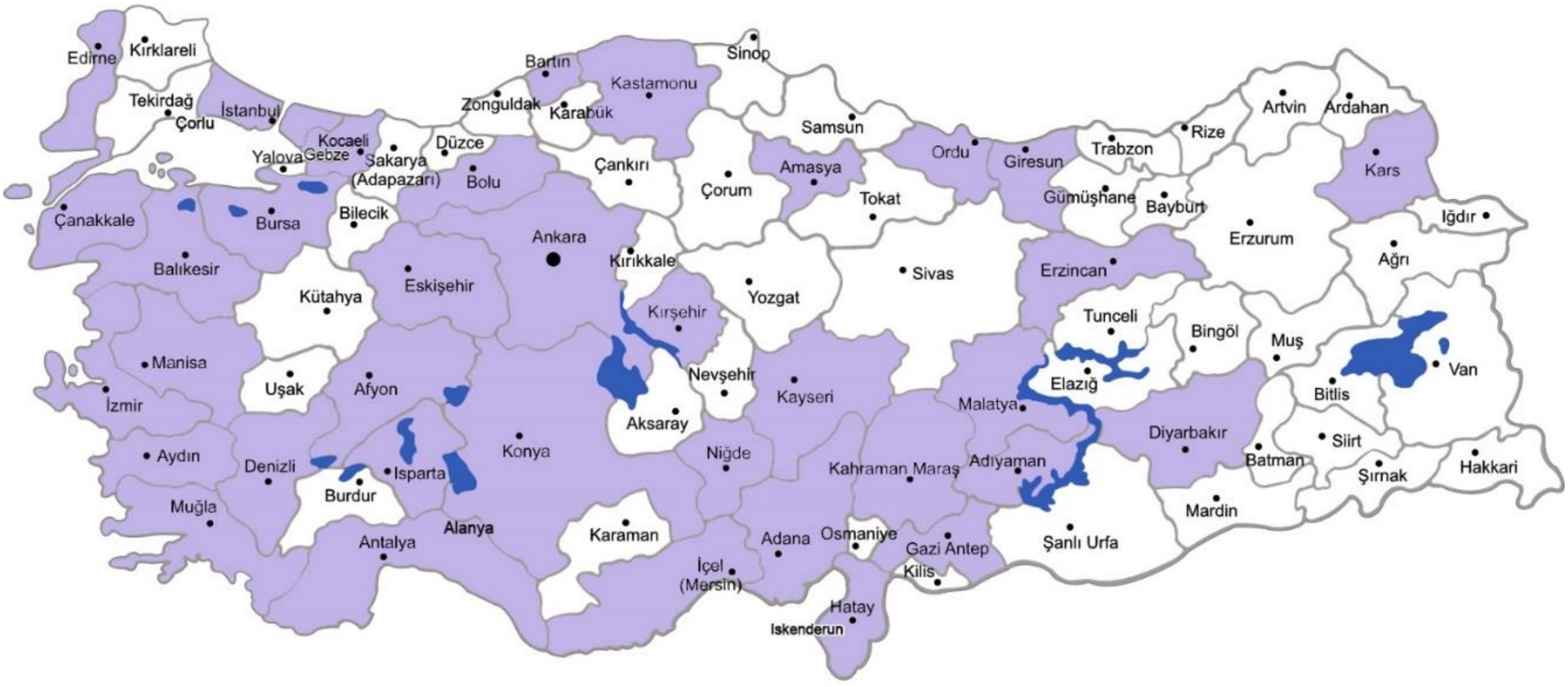

Fig. C.1 Geographical distribution of 40 cities (shaded areas)

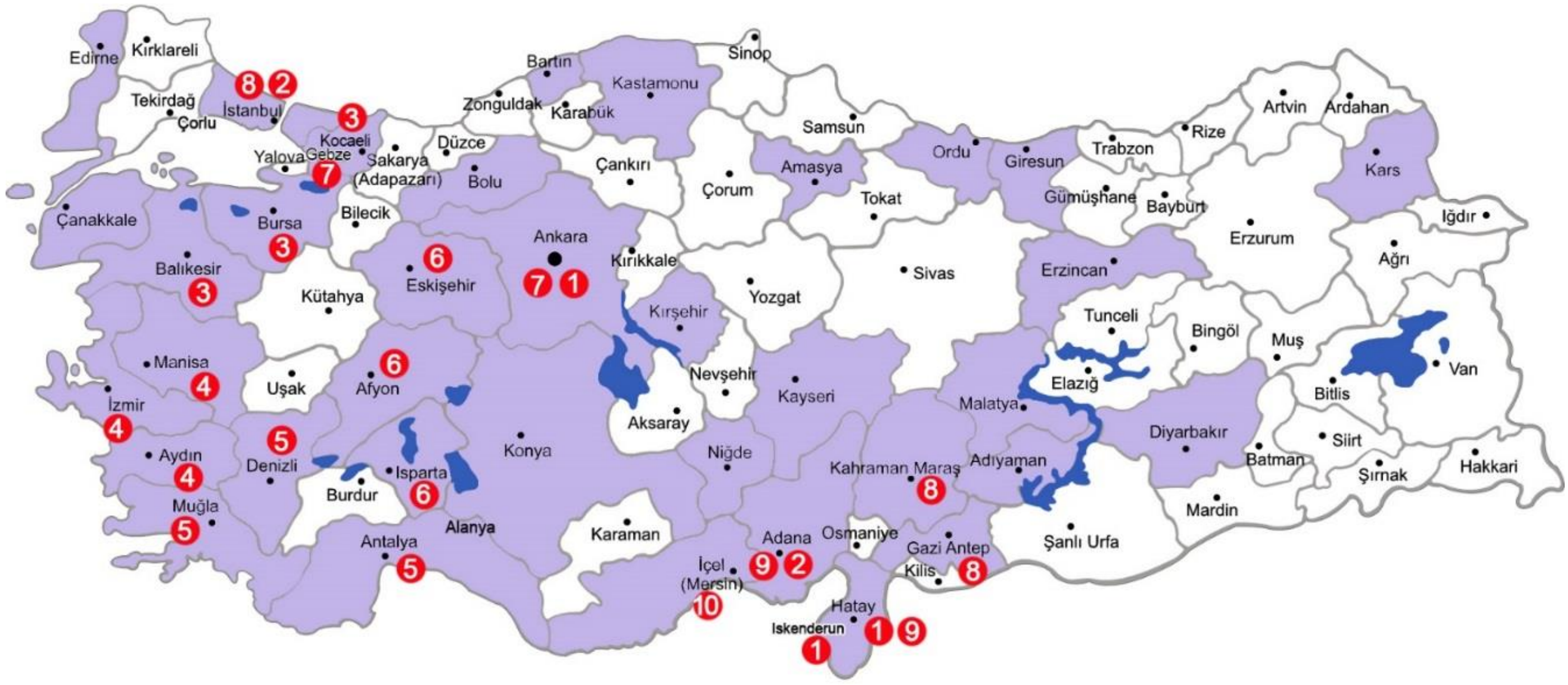

Fig. C.2 Cities with meetings